\documentclass{amsart}
\usepackage{amssymb,stmaryrd}
\usepackage{amsfonts}
\usepackage{amstext}
\usepackage{geometry}                
\usepackage{graphicx}    
\usepackage{amssymb,amsmath,epsfig}
\usepackage{scalefnt}
\usepackage{MnSymbol}
\usepackage[all]{xypic}
\usepackage{slashed}
\usepackage{extarrows}
\usepackage{stmaryrd}
\usepackage{lettrine}
\usepackage{booktabs}
\usepackage{paralist}
\usepackage{subfig}
\usepackage{color}

\newtheorem{lemma}{Lemma}[section] 
\newtheorem{proposition}[lemma]{Proposition}

\newtheorem{example}[lemma]{Example}
\newtheorem{theorem}[lemma]{Theorem}

\newtheorem{remark}[lemma]{Remark}

\renewcommand{\imath}{\mathrm{i}}

\newcommand{\CC}{\hbox{{$\mathcal C$}}}

\newcommand{\C}{\mathbb{C}}
\newcommand{\R}{\mathbb{R}}

\newcommand{\Z}{\mathbb{Z}}

\newcommand{\del}{\partial}

\newcommand{\ev}{\mathrm{ev}}
\newcommand{\coev}{\mathrm{coev}}

\newcommand{\extd}{\mathrm{d}}

\newcommand{\tens}{\mathop{{\otimes}}}

\newcommand{\id}{\mathrm{id}}

\newcommand{\<}{\langle}
\renewcommand{\>}{\rangle}

\newcommand{\cX}{\mathfrak{X}}

\allowdisplaybreaks
\begin{document}

\title{Quantum geodesic flows on graphs}

\subjclass[2000]{ }

\author{Edwin Beggs}\author{Shahn Majid}
\address{Department of Mathematics, Bay Campus, Swansea University, SA1 8EN, UK}
\address{Queen Mary University of London, 
School of Mathematical Sciences, 327 Mile End Rd, London E1 4NS, UK}
\email{e.j.beggs@swansea.ac.uk,   s.majid@qmul.ac.uk}

\maketitle

\begin{abstract}We revisit the construction of quantum Riemannian geometries on graphs starting from a hermitian metric compatible connection, which always exists. We use this method to find quantum Levi-Civita connections on the $n$-leg star graph for $n=2,3,4$ and find the same phenomenon as recently found for the $A_n$ Dynkin graph that the metric length for each outbound arrow has to exceed the length in the other direction by a multiple, here $\sqrt{n}$. We then study  quantum geodesics on graphs and construct these on the 4-leg graph and on the integer lattice line $\Z$ with a general edge-symmetric metric.\end{abstract}

\section{Introduction}

In recent works\cite{Beg:geo,BegMa:geo,BegMa:cur,LiuMa,BegMa:Z}, we have introduced a radically new way of thinking about classical geodesics which then extends to noncommutative or `quantum' Riemannian geometry (QRG) where coordinate algebras $A$ and their differential forms do not (graded) commute in the usual way. We use the formalism of the recent text \cite{BegMa}. Such QRGs could arise as a better model of spacetime due to quantum gravity effects, but what is relevant to us here is that they also arise for any graph\cite{Ma:gra}. Here, the algebra of functions $A$ is still commutative, namely functions on the vertices, but the space of differential 1-forms $\Omega^1$ is intrinsically noncommutative as an $A$-bimodule, namely spanned as a vector space by the arrows of the graph. Arrows are bilocal objects with distinct left and right actions of $A$ according to the value of a function at the source or target of the arrow. A metric in our formalism needs a bidirected graph and is then literally an assignment of a nonzero `square length' to each arrow, with a symmetry requirement. The usual choice for the latter is `edge symmetric'\cite{Ma:sq} where opposite arrows have the same length so that a metric is just an assignment to the undirected edge, i.e. again what we might intuitively think of  a metric for a graph. In this way,  once we allow noncommutativity differentials then discrete geometry sits very naturally as a special case of noncommutative differential geometry. 

What is important here is that this approach to discrete geometry is not ad-hoc but within a general formalism that applies to any $A$, with $C^\infty(M)$ at one extreme and graph geometry at another extreme. This can also be achieved within Connes formalism of noncommutative geometry e.g. via spectral triples\cite{Con}, but here we do not consider Dirac operators. Rather, the QRG formalism as in \cite{BegMa} and related works is more constructive with direct access to geometric quantities such as the quantum metric tensor $\mathfrak{g}\in \Omega^1\tens_A\Omega^1$ and a specific notion of quantum Levi-Civita connection  $\nabla:\Omega^1\to \Omega^1\tens_A\Omega^1$.  The two approaches are not incompatible and it would be an interesting question as to how to formulate comparable quantum geodesics in Connes approach. It is also fair to say that the QRG approach grew out of experience with quantum groups and a wealth of literature driven by mathematical physics as much as from purely mathematical considerations. 

Section~\ref{secpreqrg} recaps the QRG formalism in general, followed by the construction of quantum geodesics at that level in  Section~\ref{secgeo}. We then turn in  Section~\ref{secgra} to the QRG formalism in the graph case. The first order of business is that in general a quantum Levi-Civita connection (QLC) need not exist for a given metric (and if it does it need not be unique). One approach to this is to restrict the context, but another, which works in some generality, is the line in \cite{BegMa:rie}\cite[Chap.~8.3]{BegMa} to work more generally with hermitian metrics, where there are always plenty of hermitian-metric compatible connections.  A Hermitian metric and the usual notion of a metric in \cite{BegMa} are related by $*$ and a $*$-preserving hermitian-metric compatible connection is necessarily metric compatible in the usual way. We see in the graph case how this plays out and use this method to fully solve the problem for an $n$-leg star graph, finding QLC's for $n\le 4$. For $n=2$ the moduli of connections has a free phase parameter. For $n=3$, there is a unique choice of a phase or its conjugate (two connections) and for $n=4$ there is a unique connection corresponding to phase -1. We also have the same feature as in \cite{ArgMa3} for the Dynkin $A_n$ graph $\bullet$--$\bullet\cdots\bullet$--$\bullet$ that the metric lengths $g_{x\to y}, g_{y\to x}$ for each edge cannot be chosen independently but in a certain fixed ratio, in our case a factor of $\sqrt{n}$ in the outgoing direction. We then turn in Section~\ref{secgrageo} to how quantum geodesics look in the graph case and in Section~\ref{secstar} specifically for the $n$-leg star case.  

In Section~\ref{seccayley}, we revisit the particularly nice class of Cayley graphs, where the vertices form a group $G$ and arrows are given by right translation from among a set of generators $\CC$. We show that every bimodule connection in this case is of a certain `inner' form\cite{Ma:gra}\cite{BegMa} and then look at quantum geodesics on Cayley graphs. We end in Section~\ref{secZ} with the important example of $G=\Z$ for a generic edge-symmetric metric. The case of geodesics for a constant flat metric on $\Z$ was recently studied in\cite{BegMa:Z} and on $\Z_3$ in the original paper \cite{Beg:geo}.

Finally, without going into details, the idea of  quantum geodesics, of interest even in classical geometry, consists as a new way of thinking about geodesics  as follows. First, we do not  think about one geodesic at a time but rather a flow of geodesics much like in fluid dynamics, where each particle moves along a geodesic. The tangent vectors to these geodesics form a geodesic-time dependent vector field $X(t)$,  and it turns out that these obey a simple {\em geodesic velocity equation} which (classically) is just
\begin{equation}\label{claveleq} \dot X+\nabla_XX=0.\end{equation}
This says that the convective derivative of $X$ is zero in fluid terms. Next, instead of an actual (evolving) fluid particle density $\rho(t)$, we have  an evolving wave function $\psi(t)$ with $|\psi(t)|^2=\rho(t)$  thought of a  density, and we solve for $\psi$ by the {\em amplitude flow} equation
which (classically) is just
\begin{equation}\label{claamp} \dot \psi+ \kappa \psi + X(\extd \psi)=0,\quad \kappa=\tfrac12{\rm div}(X).\end{equation}
If one considers bump functions, then classically this reproduces a bump travelling with velocity $X(t)$ evaluated at the bump, i.e. a classical geodesic as expected. If $\psi$ is real-valued then this is not really different from working with $\rho$, i.e. a fluid flow language. But when $\psi$ is complex-valued then there are possible interference effects as normally associated with quantum mechanics. In the general noncommutative case $\psi(t)\in A$ is an element of a $*$-algebra and $|\psi(t)|^2$ is replaced by the positive operator $\psi(t)^*\psi(t)\ge 0$ but nevertheless behaves like a probability density when used to compute expectation values with respect to a reference positive linear functional $\int:A\to \C$. Classically, this could be the Riemannian measure of integration via weight $\sqrt{|{\rm det}(g)|}$. In our graph case, we similarly take $\int f=\sum_{x}\mu_x f(x)$ where we sum over the vertices with weights $\mu$. Quantum geodesics are defined for any measure but a natural choice is to be compatible with the quantum Riemannian geometry in the sense\cite{BegMa:cur} that $\int$ of a total divergence vanishes as would be the case classically for the Riemannian measure. Either way, we impose a reality condition on the geodesic velocity field that ensures that the flow is unitary in the sense of preserving $\int \psi^*\psi$, which we call an `improved auxiliary condition'\cite{BegMa:cur}. One of the new features of the present work is to further generalise this set-up by allowing an imaginary driving force $F$ on the right hand side of the quantum version of (\ref{claveleq}) which can alternatively be used to maintain unitarity of the amplitude flow (our previous auxiliary condition now appears $F=0$). Such auxiliary conditions would be automatic in the classical case as everything remains real without the need for any driving force. Finally, there is an ideological leap even in the classical case, namely $A$ could be or play the role of functions on {\em spacetime} not `space'. In this case $\psi$ are amplitude densities for events in spacetime not for locations of a particle in space and the time parameter, which we denoted by $t$ above, is an external {\em geodesic time}. Henceforth, we will denote the geodesic time parameter by $s$ to avoid any potential confusion with spacetime. In all cases, the `quantum mechanics-like' view of $\psi$ is interpreted with respect to geodesic time.

 \section{Preliminaries: $*$-algebraic formalism}\label{secpre}

We recap the formalism of quantum Riemannian geometry from \cite{BegMa} and then of quantum geodesics from \cite{Beg:geo,BegMa:geo, BegMa:cur}. We will be interested in the graph case but it is important that our constructions are not ad-hoc, merely a specialisation of a general framework.

\subsection{Recap of QRG formalism on an algebra}\label{secpreqrg}

A differential structure on a unital algebra $A$ is an $A$-bimodule $\Omega^1$ of differential forms equipped with a map $\extd:A\to \Omega^1$ obeying the Leibniz rule 
\[ \extd(ab)=a\extd b+ (\extd a)b\]
and such that $\Omega^1$ is spanned by elements of the form $a\extd b$ for $a,b\in A$. This can always be extended to a full differential graded `exterior algebra' though not uniquely (but with a unique maximal one). In the $*$-algebra setting, we say that we have a $*$-differential structure if $*$ extends to $\Omega$ (or at least $\Omega^1$) as a graded-involution (i.e., with an extra minus sign on swapping odd degrees) and commutes with $\extd$. 

A full formalism of quantum Riemannian geometry in this setting can be found in \cite{BegMa}. In particular, a metric means for us is an element $\mathfrak{g}\in\Omega^1\tens_A\Omega^1$ which is strongly invertible in the sense of a bimodule map $(\ ,\ ):\Omega^1\tens_A\Omega^1\to A$ obeying the usual requirements as inverse to $\mathfrak g$. This assumption forces $\mathfrak g$ in fact to be central and can be relaxed, for example by not requiring $(\ ,\  )$ to descend to the tensor product over $A$. Next, working with $\mathfrak g$,  a QLC or {\em quantum Levi-Civita connection} is a bimodule connection $(\nabla,\sigma)$ on $\Omega^1$ which is metric compatible and torsion free in the sense in the sense 
\[ \nabla \mathfrak{g}:=(\nabla\tens\id+(\sigma\tens\id)(\id\tens\nabla))\mathfrak g=0,\quad   T_\nabla:= \wedge\nabla-\extd=0\]
for a left bimodule connection $\nabla:\Omega^1\to \Omega^1\tens_A\Omega^1$ characterised by
\[ \nabla(a.\omega)=\extd a\tens \omega+ a.\nabla\omega,\quad \nabla(\omega.a)=\sigma(\omega\tens\extd a)+(\nabla\omega).a\]
where the `generalised braiding' $\sigma:\Omega^1\tens_A\Omega^1\to \Omega^1\tens_A\Omega^1$ is assumed to exist and is uniquely determined by the second equation. We also require $\mathfrak{g}$ to  be `real' and for QLC  $\nabla$ to be $*$-preserving, i.e.,
\[ \mathfrak{g}^\dagger=\mathfrak{g},\quad \sigma\,\dag\,\nabla(\xi^*)=\nabla\xi,\quad  \dagger:={\rm flip}(*\tens *).\]
Note that torsion free and $*$-preserving imply respectively for $\sigma$ that
\[ \wedge(\id+\sigma)=0,\quad \sigma\,\dag\,\sigma\,\dag=\id\]
and note that the latter implies that $\sigma$ is invertible. If these hold without necessarily the full conditions for $\nabla$ itself then we say that the connection is respectively torsion-compatible and $*$-compatible. The former says that $T_\nabla$ is a right $A$-module map and hence a bimodule map as it is already a left $A$-module map.

There is an analogous theory of right bimodule connections with left and right swapped. In this paper, we will work as usual in \cite{BegMa} with a left connection on $\Omega^1$. We also have recourse to $\cX={}_A\hom(\Omega^1,A)$, the space of left quantum vector fields defined as left $A$-module maps  $X: \Omega^1\to A$. This has a bimodule structure 
\[ (a.X.b)(\omega)= (X(\omega.a))b\]
for all $\omega\in \Omega^1$ and $a,b\in A$, and inherits a right connection. Similarly for right vector fields $\hom_A(\Omega^1,A)$.

Finally, a Hermitian metric on an $A$-bimodule $E$ is  a bimodule map $\<\ ,\ \>:E\tens\overline{E}\to A$ where we denote by $\bar E$ the same set as $E$ but the conjugate aciton of the field $\C$ and of $A$,
\[ \lambda.\bar e=\overline{\bar\lambda.e},\quad a.e=\overline{e.a^*},\quad e.a=\overline{a^*.e}\]
for all $a\in A$, where $\bar e$ denotes $e\in E$ when viewed in $\bar E$. A left (bimodule) connection $\nabla_E:E\to \Omega^1\tens_A E$ on $E$ defines
canonically a right (bimodule) connection $ \nabla_{\bar E}:\bar E\to \bar E\tens\Omega^1$ by 
\[\nabla_{\bar E}(\bar e)={\rm flip}(*\tens\bar{\ })\nabla e,\quad \sigma_{\bar E}(\omega\tens \bar e)={\rm flip}(*\tens\bar{\ })\sigma_E(e\tens\omega^*).\]
This also works just for a left connection giving a right connection. Then $\nabla_E$ is {\em hermitian metric compatible} if 
\[ \extd\<\omega,\bar\eta\>=(\id\tens\<\ ,\ \>)(\nabla\tens\id)+(\<\ ,\ \>\tens\id)(\id\tens\bar\nabla).\]

The case of immediate interest is $E=\Omega^1$ and $\nabla_E=\nabla$. Here a hermitian-metric and a round-bracket metric are equivalent data via $(\ ,\ )=\<\ ,(\ )^*\>$. We do not assume for the moment that either side descends to $\tens_A$ but we do require non-degeneracy, which amount to the map $\omega\mapsto (\ ,\omega)$ a linear isomorphism $\mathfrak{g}_2:\Omega^1\to \cX$ . Next, we let $\Omega^1$ be f.g.p. with left bases $\{e^i\},\{e_i\}$ with associated $A$-valued projection matrix $P_{ij}=\ev(e^i\tens e_j)$ and we let
$g^{ij}=\<e^i,\overline{e^j}\>=(e^i,(e^j)^*)$ is an $A$-valued matrix $g$. The initial metric being hermitian is equivalent to $g^\dagger=g$, where $(\ )^\dagger$ denotes hermitian-conjugation of a matrix. The associated lower index matrix $\tilde g=(g_{ij})$ is uniquely determined by
\[ \tilde g P=\tilde g,\quad g_{ij}e^j=\mathfrak{g}_2(e_i)^*\]
and obeys\cite[Prop.~8.31]{BegMa}
\[ g\tilde g=P,\quad \tilde g g=P^\dagger,\quad  \tilde g^\dagger=\tilde g.\]
There is an associated $\mathfrak g=g_{ij}e^i\tens_A e^j$ relevant to the QRG as above. If $(\ ,\ )$ descends to $\tens_A$ then this is central, but we do not need to assume this. Then a left connection written as $\nabla e^i=-\Gamma^i{}_j\tens e^j$ for a 1-form valued matrix $\Gamma$,  
is  hermitian-metric compatible if and only if
\[ \extd g+ \Gamma g+g \Gamma^\dagger=0\]
 of a matrix and we solve this with
\[  \Gamma=-{1\over 2}(\extd g)\tilde g + N\tilde g,\quad  N^\dagger=-N  \]
for some 1-form valued matrix $N$ as stated. This means that hermitian-metric compatible connections always exist.

In terms of $(\ ,\ )$, hermitian metric compatibility amounts to
\[ \extd(\omega,\eta)=(\id\tens(\ ,\ ))(\nabla\omega\tens\eta)+((\ ,\ )\tens\id)(\omega\tens \dagger \nabla(\eta^*)).\]
If $\nabla$ is a $*$-preserving bimodule connection and $(\ ,\ )$ descends to $\tens_A$ such that
\begin{equation}\label{metsig} (\ ,\ )\tens\id =(\id\tens (\ ,\ ))(\sigma\tens\id)(\id\tens\sigma)\end{equation}
then this is  equivalent to the usual
\begin{equation}\label{metnabla} \extd(\omega,\eta)=(\id\tens(\ ,\ ))(\nabla\omega\tens\eta+(\sigma\tens\id)(\omega\tens\nabla\eta)),\end{equation}
which in turn is equivalent to $\nabla \mathfrak g=0$, see \cite[Lemma~8.4]{BegMa}. Here (\ref{metnabla}) expresses metric compatibility as $(\ ,\ )$ intertwining $\extd$ and the tensor product bimodule connection, which implies  (\ref{metsig}).

Thus, a natural approach to solving for a QLC is to first solve for hermitian-metric compatibility, then ask for which solutions are bimodule connections, torsion free, $*$-preserving and for which $\sigma$ obeys (\ref{metsig}). In our case of an inner calculus (in the sense that $\extd=[\theta,\ ]$ for $\theta\in \Omega^1$ with $\theta^*=-\theta$), this is not necessarily best done by solving for $N$ in the above but rather by a theorem \cite{Ma:gra}\cite[Thm~8.11]{BegMa} that  if $\Omega^1$ is inner, then a left bimodule connection has the specific form
\begin{equation}\label{nablainner} \nabla \omega=\theta\tens\omega- \sigma(\omega\tens\theta)+\alpha(\omega)\end{equation}
for some bimodule maps $\sigma:\Omega^1\tens_A\Omega^1\to \Omega^1\tens_A\Omega^1$ and $\alpha:\Omega^1\to \Omega^1\tens_A\Omega^1$. 

Finally, for a quantum metric, it is usual to impose some form of symmetry. Three notions here are:
\begin{enumerate}\item quantum symmetry $\wedge(\mathfrak{g})=0$;
\item $\sigma$-symmetry $(\ ,\ )\circ\sigma=(\ ,\ )$  (but, for the pair $(\mathfrak{g},\nabla)$);
\item edge-symmetry in the graph case (see below). 
\end{enumerate}
We do not necessarily impose any of these, i.e. we are more precisely working with nondegenerate `generalised metrics'. 

\subsection{Recap of quantum geodesics on an algebra}\label{secgeo}

We do  not repeat here the origins of the theory in the 2-category of $A$-$B$-bimodules with bimodule connections, but rather just give the resulting equations for quantum geodesics.

The data we need on a differential algebra $(A,\Omega^1,\extd)$ is first of all a left bimodule connection $\nabla:\Omega^1\to \Omega^1\tens_A\Omega^1$ which can be any bimodule connection on $\Omega^1$, but for the geometric case could  be a QLC or WQLC\cite{BegMa}  with respect to a quantum metric as in Section~\ref{secpreqrg}. We assume that $\Omega^1$ is finitely generated and projective. Then the bimodule of left vector fields $\cX$ acquires a dual right connection
\[\nabla_\cX:\cX\to \cX\tens_A\Omega^1,\quad \sigma_\cX: \Omega^1\tens_A\cX\to \cX\tens_A\Omega^1.\]
 characterised by 
\[  \extd(\ev(\omega\tens X))=(\id\tens\ev)(\nabla \omega \tens X)+(\ev\tens\id)(\omega\tens\nabla_\cX X) \]
for all $\omega\in \Omega^1$ and $X\in \cX$. In the inner case where $\nabla$ is given by (\ref{nablainner}), we can give this more explicitly as follows. 
  \begin{lemma}\label{nablachi} For an inner calculus and $\Omega^1$ f.g.p. with dual bases $\{e^i\},\{e_i\}$, the dual bimodule connection can be given explicitly
\begin{align*}
\nabla_\cX X &= - X\tens\theta+(\id\tens\id\tens\ev)(e_i\tens\sigma(e^i\tens\theta)\tens X-e_i\tens\alpha(e^i)\tens X)  \cr
\sigma_\cX(\xi\tens X) &= (\id\tens\id\tens\ev)(e_i\tens\sigma(e^i\tens\xi)\tens X)
\end{align*}
 (here $\coev(1)=e_i\tens e^i\in \cX\tens\Omega^1$, summation implicit).
 \end{lemma}

We can equally well formulate metric compatibility of $\nabla$ in terms of $\nabla_\cX$. If $\sigma$ is invertible and $\nabla$ is metric compatible then using (\ref{metsig}) we also have a `mixed' form of metric compatibility that the map $\mathfrak{g}_2$ intertwines the right bimodule connections $\nabla_\cX$ and $\sigma^{-1}\circ\nabla$,   \begin{equation}\label{sigmachiinv}
\nabla_\cX=(\mathfrak{g}_2\tens\id)\circ\sigma^{-1}\circ\nabla\circ\mathfrak{g}_2,\quad \sigma_\cX{}^{-1}=(\id\tens \mathfrak{g}_2) \circ\sigma_{\Omega^1}\circ( (\mathfrak{g}_2)^{-1}  \tens\id)
\end{equation}
 where the second half follows from the first half as a morphism of bimodule connections. 
 
Next, for any right connection $\nabla_\cX$ with $\sigma_\cX$ is invertible,  $\hat\nabla=\sigma_\cX^{-1}\nabla_\cX$ with $\hat\sigma=\sigma_\cX^{-1}$ is  a left bimodule connection on $\cX$, 
which,  allows us to define a divergence of $X\in\cX$ as
\[ {\rm div}:=\ev\circ\hat\nabla=\ev\circ\sigma_\cX^{-1}\nabla_\cX,\quad \ev:\Omega^1\tens_A\cX\to A,\quad \ev(\omega\tens_A X)=X(\omega)\]
 for all $\omega\in \Omega^1, X\in \cX$. Here $\ev$ is a bimodule map and one can check that 
 \[ {\rm div}(aX)=X(\extd a)+ a{\rm div}(X).\]
Moreover, if $\nabla$ is metric compatible with $\sigma$ invertible then we can use (\ref{sigmachiinv}) to find
\begin{equation}\label{div2} {\rm div}(\mathfrak{g}_2(\omega))=\ev\circ\sigma_\cX^{-1}\circ\nabla_\cX\circ\mathfrak{g}_2(\omega)= \ev\circ(\id\tens \mathfrak{g}_2)\circ\sigma\circ \sigma^{-1}\circ\nabla\circ\mathfrak{g}_2(\omega) =(\ ,\ )\nabla \omega\end{equation}
on substituting both halves of (\ref{sigmachiinv}). So in the metric compatible case, the divergence of a vector field matches the natural codifferential of the corresponding 1-form. 

We are now ready for quantum geodesics. Although the theory is more general, we will focus on `wave functions' $\psi\in E=C^\infty(\R, A)$, where $s\in\R$ will be the `geodesic time' parameter. Likewise, we let $X_s$ be a time-dependent 1-form om $A$ and $\kappa_s$ another time-dependent element of $A$. These obey the {\em geodesic velocity equations} if
\begin{equation}\label{veleqX} \dot X_s +[X_s,\kappa_s]+(\id\tens X_s)\nabla_\cX(X_s)=0\end{equation}
where dot means differential with respect to $s$. 

Next, we require $\int: A\to \C$ to be a non-degenerate positive linear functional (so we can think of it as a probability measure if normalisable so that $\int 1=1$) and define  ${\rm div}_{\int}(X)$ of a vector field by
\[ \int a\,  {\rm div}_{\int}(X)+ \int X(\extd a)=0\]
for all $a\in A$. If this is the same as the geometric divergence then we say that $\int$ is {\em divergence compatible} (with $\nabla$), which is equivalent to 
\[ \int{\rm div} X=0\]
for all $X\in \cX$. We do not necessarily assume this, however, as it can be quite restrictive in noncommutative geometry. We do require that the geodesic velocity field and $\kappa$ at each $s$ to obey the {\em unitarity conditions}  
 \begin{equation}\label{unitarity} \int \kappa^* a+a\kappa+X(\extd a)=0,\quad \int X(\omega^*)=X(\omega)^*\end{equation}
for all $a\in A$ and all $\omega\in\Omega^1$. Note that if the second equation implies (one says that $X$ is {\em real with respect to $\int$}) then we can canonically solve the first equation in the pair by $\kappa={1\over 2}{\rm div}_{\int}(X)$, see after \cite[Def.~4.5]{BegMa:cur}. It is not automatic that if $X_s$ is initially real with respect to $\int$ that it necessarily remains so under the geodesic velocity equation, we have to impose this as a further `{\em improved auxiliary equation}' obtained as the difference between (\ref{veleqX}) and its conjugate under the reality assumption. It may be that this is not possible, in which case one can always add a time-dependent driving force $F\in \cX$ to (\ref{veleqX}) for which a natural choice is to take $\imath F$ real with respect to $\int$ also. This then uniquely defines $F$ as an external force needed to maintain $X$ real with respect $\int$ as it evolves. The improved auxiliary equation appears in this extended framework as $F=0$. We will see how this greater freedom plays out in the graph case.

Given such a geodesic velocity vector field, we then require the {\em amplitude flow equation} 
\[ \dot\psi= -\psi\kappa_s-X_s(\extd \psi),\]
where $\extd$ acts on $\psi(s)\in A$ and dot is with respect to $s$ as before. The above conditions ensure that $\int\psi^*\psi$ is constant in time, which is needed for a probabilistic interpretation. 

We also assume that $\Omega^1$ is a $*$-calculus and $\nabla$ is $*$-preserving. In the nicest case, there is also a 
 $*$ operation on $\cX$, namely\cite{BegMa:cur} when $\int$ is a (twisted) trace in the sense 
\[ \int ab=\int\varsigma(b)a\]
for all $a,b\in A$ with respect to an algebra automorphism $\varsigma$ that extends to a map $\varsigma:\Omega^1\to\Omega^1$ by $\varsigma(a.\extd b)=\varsigma(a).\extd \varsigma(b)$.

\begin{theorem} \label{thmXstar}\cite{BegMa:cur} If a positive linear functional $\int$ is a twisted trace with respect to $\varsigma$ then $*$ on $\cX$ defined by
\[
 X^*(\omega)=\big(\ev\circ\sigma_\cX^{-1}(X\tens\omega^*)\big)^*
\]
for all $\omega\in\Omega^1$ obeys 
\[ X^{**}=X,\quad (a\,X)^* =X^*a^*,\quad (Xa)^* =a^*\, X^*,\quad \int(X^* (\omega^*))=\int(X(\varsigma(\omega)))^*\]
for all $a\in A, X\in \cX$. If in addition $\int$ is divergence compatible,  $X$ is `real' in the sense $X^*=\varsigma\circ X\circ\varsigma^{-1}$ and 
\[ \kappa={1\over 2}{\rm div}(X)\]
as before then the unitarity conditions (\ref{unitarity})  hold.
\end{theorem}
One also has 
\[ {\rm div}(X^*)=*\circ{\rm div}(X)=\varsigma{\rm div}(X)\]
in the case of real $X$  in the  sense stated. We therefore impose this reality condition on $X$ (which for a trace just means $X^*=X$). It amounts to the same as imposing reality of $X$ with respect to $\int$ but now with a geometric meaning behind that as self-adjoint with respect to an involution.  This nice setting is, however, quite restrictive and does not always apply. On solving the geodesic velocity equations,  we again want an initially real $X_s$ to remain real and can impose this as an auxiliary condition or more generally ensure it automatically by a generated imaginary driving force $F$. 

\section{QRG formalism applied to graphs}\label{secgra}

Here, we develop the quantum Riemannian geometry of graphs, building on \cite{Ma:gra} but now coming at it from the more general framework of Hermitian-metric compatible connections as in \cite{BegMa:rie}. If $A=C(X)$ is the algebra of all complex-valued functions on a discrete set $X$ then its possible $\Omega^1$ are known to be in 1-1 correspondence with directed graphs with vertex set $X$, where we exclude self-arrows and multiple arrows.
 So a graph is just a discrete set with differential structure of order 1. We shall abuse notation by using $x\to y$ both as an arrow in the graph and as a logical statement `there is an arrow from $x$ to $y$'.  The associated calculus has 1-forms spanned as a vector space by the arrows, or equally with basis $\{\omega_{x\to y}\}$ labelled by arrows. The bimodule structure and exterior derivative are
\[ f.\omega_{x\to y}=f(x)\omega_{x\to y},\quad
\omega_{x\to y}.f=f(y)\omega_{x\to y},\quad  \extd f=\sum_{x\to y}(f(y)-f(x))\omega_{x\to y}.\]
 This forms a $*$-calculus with $f^*(x)=\overline{f(x)}$ and $\omega_{x\to y}^*=-\omega_{y\to x}$ and is inner with $\theta=\sum_{x\to y}\omega_{x\to y}$. We will also need higher degree forms and the most natural one that is functorial (i.e. defined for all graphs) is 
 \[ \Omega_{min}:\quad \forall x,z\in X:\quad \sum_{y: x\to y \to z} \omega_{x\to y}\wedge\omega_{y\to z}=0\]
 where we add the quadratic relations shown. This is explained in \cite[Chap 1]{BegMa} as a natural quotient of the maximal prolongation that
 remains inner by  $\theta$ in all degrees.

A Hermitian metric in this context and in the simplest case where it descends to $\tens_A$, is given by coefficients $\lambda_{x\to y}\in \R\setminus \{0\}$ on the arrows and has the form
\begin{align} \label{prus1}
\< \omega_{x\to y},\overline{ \omega_{z\to w}  } \> = \lambda_{x\to y}\,\delta_{x,z}\,\delta_{y,w}\, \delta_x, 
\end{align}
with the natural choice being $\lambda_{x\to y}>0$. 
Here the $\delta_{x,z}$ is given by the inner product being a bimodule map, and the $\delta_{y,w}$ comes from the assumption of descent to $\tens_A$.
This corresponds via $*$ to a usual metric  
\[\mathfrak{g}=-\sum_{x\to y} g_{x\to y}\omega_{x\to y}\tens\omega_{y\to x},\quad (\omega_{x\to y},\omega_{y'\to x'})=-\lambda_{x
\to y} \delta_{x,x'}\delta_{y,y'}\delta_x,\quad  \lambda_{x\to y}={1\over g_{y\to x}}\]
where the natural choice with $\lambda_{x\to y}>0$ is for metric coefficients here to be negative. This is a little counter-intuitive but fits with $\omega_{x\to y}^*=-\omega_{y\to x}$. Note that for constructing a QLC, we do not care about the overall normalisation of the metric so from the non-hermitian point of view we would normally omit this minus sign.

Next, since the graph calculus as well as $\Omega_{min}$ are inner, we can use the form (\ref{nablainner}) from \cite{Ma:gra}\cite[Thm.~8.11]{BegMa} for any bimodule maps $\sigma:\Omega^1\tens_A
\Omega^1\to \Omega^1\tens_A\Omega^1$ and $\alpha:\Omega^1\to \Omega^1\tens_A\Omega^1$. To this we add the idea that Hermitian-metric compatible ones always exist\cite[Chapter~8.3]{BegMa} so it is easier to solve for these first and then require $*$-preserving and torsion-freeness. We say that $\nabla$ is {\em inner} if $\alpha=0$.

\begin{figure}
\unitlength 0.3 mm
\begin{picture}(110,40)(0,20)
\linethickness{0.3mm}
\multiput(20,30)(0.12,0.16){125}{\line(0,1){0.16}}
\linethickness{0.3mm}
\multiput(35,50)(0.12,-0.16){125}{\line(0,-1){0.16}}
\linethickness{0.3mm}
\put(20,30){\line(1,0){30}}
\linethickness{0.3mm}
\put(90,30){\line(0,1){20}}
\linethickness{0.3mm}
\put(90,30){\line(1,0){20}}
\linethickness{0.3mm}
\put(110,30){\line(0,1){20}}
\linethickness{0.3mm}
\put(90,50){\line(1,0){20}}
\put(14,30){\makebox(0,0)[cc]{$x$}}
\put(57,29){\makebox(0,0)[cc]{$y$}}
\put(84,30){\makebox(0,0)[cc]{$x$}}
\put(117,29){\makebox(0,0)[cc]{$y$}}
\put(35,55){\makebox(0,0)[cc]{$z$}}
\put(85,51){\makebox(0,0)[cc]{$z$}}
\put(115,51){\makebox(0,0)[cc]{$u$}}
\end{picture}
\caption{configurations in the undirected graph defining $N$ and $L$ in Prop.\, \ref{propHM}}\label{yu9}
\end{figure}
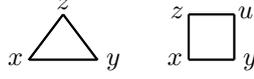

\begin{proposition} \label{propHM}
For nondegenerate metrics (i.e.\ all $\lambda_{x\to y}\neq 0$) there is a 1-1 correspondence between Hermitian-metric compatible left bimodule connections on $\Omega^1$ and the following data:

\smallskip
\noindent (1) $N_{   x\to y  , z\to y  }\in\C$ for all triangles in Fig.\, \ref{yu9} obeying $N_{   x\to y  , z\to y  }{}^*=N_{ z\to y ,  x\to y  }$.
\newline (2) $L_{   x\to y  , z\to u  }\in\C$ for all squares in Fig.\, \ref{yu9} obeying $L_{   x\to y  , z\to u  }{}^*=L_{  z\to u, x\to y    }$. Note that the squares are allowed to collapse diagonally opposite vertices, i.e.\ we can have $z=y$ or $x=u$.

\smallskip \noindent
The connection is then defined by (\ref{nablainner}) using
\begin{align} \label{jutr3}
\alpha(\omega_{x\to y}) &= \sum_{ z \, : \,  x\to z\to y  }  \lambda_{x\to y }\,N_{   x\to y  , z\to y  }\,  \omega_{x\to z} \tens \omega_{z\to y} \cr
\sigma(\omega_{x\to y}\tens\omega_{y\to u}) &=   \sum_{ z \, : \,  x\to z\to u  } \lambda_{x\to y }\,L_{   x\to y  , z\to u  }\,  \omega_{x\to z} \tens \omega_{z\to u}\ .
\end{align}
\end{proposition}
\proof  We set $K=\alpha-\sigma_\theta$, where $\sigma_\theta(\omega)=\sigma(\omega\tens\theta)$, and a general left connection on $\Omega^1$ can be written
\[
\nabla(\omega_{x\to y}) =\theta\tens \omega_{x\to y} + K(\omega_{x\to y})
\]
where we write the left module map $K$ in terms of the basis as
\begin{align} \label{mappdeff}
K(\omega_{x\to y}) = \sum_{ z \, : \,  x\to z\to y  } K_{   x\to y  , z\to y  }\,  \omega_{x\to z} \tens \omega_{z\to y}
+  \sum_{ z,u \, : \,   x\to z\to u \leftarrow y  } K_{   x\to y  , z\to u  }\,  \omega_{x\to z} \tens \omega_{z\to u}
\end{align}
and we take the . We can split this as
\begin{align} \label{jutr}
\alpha(\omega_{x\to y}) &= \sum_{ z \, : \,  x\to z\to y  } K_{   x\to y  , z\to y  }\,  \omega_{x\to z} \tens \omega_{z\to y} \cr
\sigma(\omega_{x\to y}\tens\omega_{y\to u}) &=  - \sum_{ z \, : \,  x\to z\to u  } K_{   x\to y  , z\to u  }\,  \omega_{x\to z} \tens \omega_{z\to u}\ .
\end{align}
For a given basis order of the directed arrows $x\to y$ we can make $K$ into a matrix
(with entries zero outside the ranges in (\ref{mappdeff}))
 and also let $G$ be the Hermitian diagonal matrix with entries $\lambda_{x\to y}\in\R$.
Now Hermitian metric compatibility is given by the matrix $KG$ being Hermitian. If $G$ is invertible (i.e.\ all $\lambda_{x\to y}\neq 0$) then all solutions of this are of the form $K=GL$ for $L$ a Hermitian matrix. \endproof

We use the relations for $\Omega^2_{min}$. Note that if we were to add more relations into the calculus then the conditions for a torsion free connection will typically be weakened.

\begin{proposition}  \label{propHR}
Hermitian metric compatible connections in Proposition~\ref{propHM} are torsion compatible precisely when 
\begin{align} \label{nondeg9}
g_{y\to x} + g_{u\to y} = g_{u\to z} + g_{y\to u}
\end{align}
holds for every nondegenerate square in the graph (see Figure~\ref{yu9}). In this case the connections are given in terms of an assignment
$Q_{x\to y}\in\C$ for every directed edge with
$Q_{x\to y}{}^*=Q_{y\to x}$ and 
\[
L_{   x\to y  , z\to u  } =Q_{u\to y}+(1-\delta_{u,x} )g_{u\to y}+(1-\delta_{y,z} )g_{y\to x}.
\]
In addition, the torsion for $\Omega_{min}$ vanishes precisely when $\alpha$ is of the form
\[
\alpha(\omega_{x\to y}) = \lambda_{x\to y}\, b_y\, \sum_z \omega_{x\to z}\tens\omega_{z\to y}
\]
for some real function $b_y$ on the vertices, i.e.\ $N_{   x\to y  , z\to y  }=b_y$. 
\end{proposition}
\proof First write
\[
\nabla\omega_{x\to y}=\theta\tens\omega_{x\to y}+\omega_{x\to y}\tens\theta +\alpha(\omega_{x\to y})-(\sigma+\id)(\omega_{x\to y}\tens\theta)
\]
We require that the following is in the kernel of $\wedge$ 
\[
(\sigma+\id)(\omega_{x\to y}\tens\omega_{y\to u})=
 \sum_{ z \, : \,  x\to z\to u  } ( \lambda_{x\to y }\,L_{   x\to y  , z\to u  } +\delta_{y,z})\,  \omega_{x\to z} \tens \omega_{z\to u}
\]
so the bracket must be independent of $z$, so we set
\[
\lambda_{x\to y }\,L_{   x\to y  , z\to u  } +\delta_{y,z}=\lambda_{x\to y }\,M_{xyu}\ .
\]
Now the Hermitian condition on the $L$s gives for every square
\begin{align} \label{nondeg}
M_{xyu}-g_{y\to x}\delta_{z,y}=M_{zuy}^*-g_{u\to z}\delta_{x,u}\ .
\end{align}
On totally degenerate squares (a single edge given by putting $z=y$ and $u=x$) we get
\begin{align}
M_{xyx}-g_{y\to x}=M_{yxy}^*-g_{x\to y}
\end{align}
and to solve this we put $M_{xyx}=g_{y\to x} +Q_{x\to y}$ where $Q_{x\to y}{}^*=Q_{y\to x}$. Next we look at the singly degenerate square where $z=y$ and $u\neq x$
\begin{align}
M_{xyu}=M_{yuy}^*+g_{y\to x}= Q_{u\to y}+g_{u\to y}+g_{y\to x}
\end{align}
and at the singly degenerate square where $z\neq y$ and $u=x$
\begin{align}
M_{zxy}= M_{xyx}^*   +   g_{x\to z} = Q_{y\to x} +   g_{y\to x} + g_{x\to z}
\end{align}
and on changing the letters these are seen to give the same information. Now we have the nondegenerate squares which give on substitution of our previous results into (\ref{nondeg}) gives $M_{xyu}=M_{zuy}^*$ and this gives (\ref{nondeg9}). 
Next we check that the resulting $L$ is Hermitian, i.e.
\[
(1-\delta_{u,x} )g_{u\to y}+(1-\delta_{y,z} )g_{y\to x}=(1-\delta_{y,z} )g_{y\to u}+(1-\delta_{u,x} )g_{u\to z}
\]
which we reorder as
\[
(1-\delta_{u,x} )(g_{u\to y}  - g_{u\to z})=(1-\delta_{y,z} )(g_{y\to u}
-g_{y\to x})
\]
which gives no new information. From the formula for $\alpha$ in (\ref{jutr3}) we see that torsion free implies that $N_{   x\to y  , z\to y  }$ is independent of $z$. Then the Hermitian condition on $N$ shows that it is also independent of $x$. \endproof

\begin{figure}[h] 
\unitlength 0.4 mm
\begin{picture}(75,40)(0,17)
\linethickness{0.3mm}
\put(30,20){\line(0,1){40}}
\linethickness{0.3mm}
\put(30,60){\line(1,0){40}}
\linethickness{0.3mm}
\put(70,20){\line(0,1){40}}
\linethickness{0.3mm}
\put(30,20){\line(1,0){40}}
\linethickness{0.3mm}
\multiput(30,20)(0.12,0.12){333}{\line(1,0){0.12}}
\put(50,40){\makebox(0,0)[cc]{$\bullet$}}

\put(75,60){\makebox(0,0)[cc]{$x$}}
\put(25,20){\makebox(0,0)[cc]{$u$}}
\put(43,41){\makebox(0,0)[cc]{$z$}}
\put(25,60){\makebox(0,0)[cc]{$y$}}
\put(75,20){\makebox(0,0)[cc]{$v$}}
\put(20,45){\makebox(0,0)[cc]{}}
\end{picture}
\caption{configurations in the undirected graph for Prop.\, \ref{propHvvy}}\label{yu77}
\end{figure}
\begin{figure}[h] 
\unitlength 0.4 mm
\begin{picture}(75,40)(0,17)
\linethickness{0.3mm}
\put(30,20){\line(0,1){40}}
\linethickness{0.3mm}
\put(30,60){\line(1,0){40}}
\linethickness{0.3mm}
\put(70,20){\line(0,1){40}}
\linethickness{0.3mm}
\put(30,20){\line(1,0){40}}
\linethickness{0.3mm}
\multiput(30,20)(0.12,0.12){333}{\line(1,0){0.12}}

\put(75,60){\makebox(0,0)[cc]{$x$}}
\put(25,20){\makebox(0,0)[cc]{$y$}}
\put(25,60){\makebox(0,0)[cc]{$z$}}
\put(75,20){\makebox(0,0)[cc]{$v$}}
\put(20,45){\makebox(0,0)[cc]{}}

\end{picture}
\caption{configurations in the undirected graph for Prop.\, \ref{propHvvy}}\label{yu99}
\end{figure}

\begin{proposition} \label{propHvvy}
Given a Hermitian-metric compatible connection on a graph, $*$-compatibility (i.e.\ $\sigma\,\dag\,\sigma\,\dag=\id$) reduces to
\[
\sum_z \lambda_{u\to y} \lambda_{x\to z} L_{z\to x,u\to y} L_{x\to z,v\to u}=\delta_{v,y}
\]
where Figure~\ref{yu77} shows the configurations of the points (which may degenerate). 
If this holds then the connection is $*$-preserving  if and only if 
\[
N_{y\to x,v\to x}=-\lambda_{x\to y} \sum_z N_{x\to y,z\to y}{}^* \lambda_{y\to z} L_{y\to z,v\to x}
\]
where Figure~\ref{yu99} shows the configurations of the points.
\end{proposition}
\proof  It is known \cite[Prop.~8.11]{BegMa} that if $\nabla$ is $*$ compatible then the inner part is also $*$-preserving hence for the connection to be $*$-preserving in the presence of $\alpha$ we just need $\sigma\,\dag\,\alpha(\xi^*)=\alpha(\xi)$, which reduces as shown for $\xi=\omega_{y\to x}$. \endproof

In both propositions, one can set  $\alpha=0$ so full torsion and $*$-preserving hold for these such inner connections as soon as they are torsion and $*$-compatible. This applies similarly to all inner calculi\cite[Prop.~8.11]{BegMa}. Also, when our connection is torsion free and $*$-preserving then by our general comments it is metric compatible in the usual (non hermitian sense) also, i.e. a QLC. 

\subsection{QRG of star  graphs}\label{secqrgstarAn}
\begin{figure}
\[\includegraphics{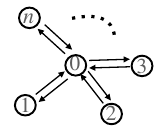}\]
\caption{\label{yu0} $n$-star graph which we decorate with metric values $g_{ i\to 0}=\lambda^{-1}_{0\to i}$ and    $g_{0 \to i}=\lambda^{-1}_{i\to 0}$  for $i\in\{1,2,\cdots n\} $  on the arrows.}
\end{figure}

We consider the $n$-star graph with $n$ vertices labelled $\{1,2,\dots,n\}$ joined to a central vertex labelled $0$ as shown 
 in Figure~\ref{yu0}. We use the notation $i,j,k\in \{1,2,\dots,n\}$ for the exterior vertices.

 \begin{theorem} \label{thmsolstar}   For the $n$-star graph, hermitian-metric compatible and $*$-compatible connections exist only if  \begin{align*}
 \frac{\lambda_{0\to i}}{   \lambda_{i\to 0} }=\frac{1}{\sqrt{n}}
\end{align*}
and are necessarily $*$-preserving with $\alpha=0$. Using $\Omega_{min}$, torsion free such connections are the moduli of QLCs and exist only for $n\le 4$. They consist of (a)  solutions
\[
\lambda_{0\to k} \, L_{0\to k,j\to 0} 
= 1-\delta_{j,k}+ {s_k^{-1}\over\sqrt{n}};\quad  s_k=\begin{cases}  e^{3\imath\pi\over 4}& n=2\\ e^{5 \imath\pi\over 6}& n=3 \\ -1 & n=4\end{cases}\]
up to complex conjugation (i.e. 2 connections for each $n$), where  $s_k$ is independent of $k$, and (b), when $n=2$,  more general solutions of the form
\[ |s_1|=1, \quad  s_2=-{s_1+\sqrt{2}\over \sqrt{2}s_1+1}. \]
Here, $s_1$ is a free phase parameter and $s_2$ is obtained from $s_1$ by a M\"obius transform that maps phases to phases and has the special values in (a) as its fixed points. Moreover,   $L_{j\to 0,0\to k}=L_{0\to k,j\to 0}{}^*$ fully specifies $\sigma$. 
\end{theorem}
\proof From Proposition~\ref{propHvv}
we have the two conditions, for $i,j,k\in\{1,2,\dots,n\}$
\begin{align} \label{huer}
\sum_k \lambda_{0\to i} \lambda_{0\to k} L_{k\to 0,0\to i} L_{0\to k,j\to 0} &=\delta_{j,i}  \cr
 \lambda_{i\to 0} \lambda_{j\to 0} L_{0\to j,i\to 0} L_{j\to 0,0\to i} &=1\ .
\end{align}
Note that we must have non-vanishing $L$'s from the second equation.
Substituting one condition into the other gives for all $i,j$
\begin{align*}
 \frac{\lambda_{0\to i}}{   \lambda_{i\to 0} } \, 
\sum_k \frac{\lambda_{0\to k}}{   \lambda_{k\to 0}   }  \ \frac{L_{0\to k,j\to 0} }
{   L_{0\to k,i\to 0}      }     &=\delta_{j,i}  
\end{align*}
In the case where $i=j$ the ratio of the $L$s cancels and it follows that the fraction 
$ \frac{\lambda_{0\to i}}{   \lambda_{i\to 0} }$ is independent of $i$, and we can solve to find the first stated condition. Now for $i\neq j$ we have the restriction
  \begin{align} \label{trub}
\sum_k   \frac{L_{0\to k,j\to 0} }
{   L_{0\to k,i\to 0}      }     &=0\ .
\end{align}
We then have $L_{j\to 0,0\to k}=L_{0\to k,j\to 0}{}^*$ and $\alpha=0$. Note that so far we have used a subset of the
metric preserving relations, but if solutions exist then by \cite[Thm~8.11]{BegMa}, due to the calculus being in the inner case, they have a form given by $\sigma$ and a bimodule map $\alpha$. The latter, however, must be zero since there are no triangles in the graph. In this case $*$-compatible implies $*$-preserving. 

Next, proceeding with this partial solution for $L$s,  torsion compatibility (which in the present inner case is equivalent to torsion freeness)  becomes
\[
L_{0\to k,j\to 0} = Q_{0\to k}+(1-\delta_{j,k})\, g_{k\to 0}
\]
and substituting this in the sum (\ref{trub}) gives for all $i\neq j$,
\begin{align} \label{vyiu}
n-2+ \frac{ Q_{0\to i}+g_{i\to 0} }{ Q_{0\to i} }+  \frac{ Q_{0\to j} }{ Q_{0\to j}+g_{j\to 0} } \ =\ 0.
\end{align}
The solutions of this simplify when $n\neq 2$ as we can show that $Q_{0\to 2}/g_{2\to 0} =Q_{0\to 1}/g_{1\to 0} $.
In this case the solutions are, using the same choice of sign $\epsilon=\pm 1$ for all $i$,
\[
\frac{ Q_{0\to i}+g_{i\to 0} }{ Q_{0\to i} } =  \frac{-n+2+ \epsilon\,\sqrt{n^2-4n }    }{2}, 
\]
which simplifies to
\[
 \frac{ g_{k\to 0} }{ Q_{0\to k} } = \frac{-n+ \epsilon\,\sqrt{n^2-4n }    }{2}.
\]
Proposition~\ref{propHR} then gives us
\[ \lambda_{0\to k} \, L_{0\to k,j\to 0} =\frac{2}{-n+ \mathrm{i}\,\epsilon\,\sqrt{4n-n^2 }    }+1-\delta_{jk},\]
which we write as stated. Now we check consistency with the remaining relations, which can be summarised by, for $i\neq j$,
\[
| \lambda_{0\to i}\, L_{0\to i,i\to 0}|=\frac{1}{\sqrt{n}}\ ,\quad 
\lambda_{i\to 0}\, \lambda_{j\to 0}\, L_{0\to j,i\to 0}\, L_{0\to i,j\to 0}{}^*=1\ .
\]
These can be checked  to hold for $n=3,4$ and fail for $n\ge 5$. 
\newline In the $n=2$ case (\ref{vyiu}) becomes
\begin{align*} 
\frac{ Q_{0\to 1}+g_{1\to 0} }{ Q_{0\to 1} }+  \frac{ Q_{0\to 2} }{ Q_{0\to 2}+g_{2\to 0} } \ =\ 0
\end{align*}
and solving this for $Q_{0\to 2}/g_{2\to 0} $ in terms of  $Q_{0\to 1}/g_{1\to 0} $ gives the solution as stated.
\endproof

For the type (a) solutions, if we write $s=s_k$ for the common value then $\sigma$ in Theorem~\ref{thmsolstar} comes out uniformly as
\begin{align*}
\sigma(\omega_{0\to i}\tens \omega_{i\to 0}) &= {s^{-1}\over\sqrt{n}}\  \omega_{0\to i}\tens \omega_{i\to 0}
+\Big(1+ {s^{-1}\over\sqrt{n}}\Big)\sum_{j\ne i}\omega_{0\to j}\tens \omega_{j\to 0}\ ,\cr
\sigma(\omega_{i\to 0}\tens \omega_{0\to i}) &=   s\  \omega_{i\to 0}\tens \omega_{0\to i}\ ,\cr
\sigma(\omega_{i\to 0}\tens \omega_{0\to j}) &=- \frac{ \lambda_{0\to i}  }{ \lambda_{0\to j}  }\  s^{-1}\  \omega_{i\to 0}\tens \omega_{0\to j}
\end{align*}
for all $j\ne i$.  The associated connection, which is necessarily of inner type governed by $\sigma$, is 
\begin{align*} \nabla\omega_{0\to i}&=\sum_j\omega_{j\to 0}\tens\omega_{0\to i}-\sigma(\omega_{0\to i}\tens\omega_{i\to 0})\\
&=\omega_{i\to 0}\tens\omega_{0\to i}- {s^{-1}\over\sqrt{n}}\omega_{0\to i}\tens\omega_{i\to 0}+\sum_{j\ne i}\Big(\omega_{j\to 0}\tens\omega_{0\to i}- (1+{s^{-1}\over\sqrt{n}})\omega_{0\to j}\tens\omega_{j\to 0}\Big)
\end{align*}
\begin{align*} \nabla\omega_{i\to 0}&=\omega_{0\to i}\tens\omega_{i\to 0}-\sum_j\sigma(\omega_{i\to 0}\tens\omega_{0\to j})=\omega_{0\to i}\tens\omega_{i\to 0}- s \omega_{i\to 0}\tens\omega_{0\to i}+\sum_{j\ne i}{\lambda_{0\to i}\over\lambda_{0\to j}}s^{-1}\omega_{i\to 0}\tens\omega_{0\to j}.
\end{align*}
 We can potentially have more $*$-preserving solutions, e.g. for larger $n$, if we quotient $\Omega_{min}$ further so that being torsion-free is less restrictive. 

\begin{remark}\rm The $n=2$ case is the same as the  $A_3$ Dynkin graph $\bullet$--$\bullet$--$\bullet$ treated in \cite{ArgMa3} and a careful comparison noting that the metric coefficients there relate to ours by
\[ h_1={1\over\lambda_{1\to 0}},\quad h_2={1\over \lambda_{0\to 2}},\]
shows that our constraints on the ratios of inbound and outbound metric coefficients as the same as found there. We also obtain exactly the  same connection $\nabla$ as in \cite{ArgMa3} with $s=s_1$ as the free modulus parameter found there.
\end{remark}

\section{Quantum geodesics on finite graphs}\label{secgrageo}

\begin{figure}
\unitlength 0.3 mm
\begin{picture}(110,40)(0,20)
\linethickness{0.3mm}
\multiput(20,30)(0.12,0.16){125}{\line(0,1){0.16}}
\linethickness{0.3mm}
\multiput(35,50)(0.12,-0.16){125}{\line(0,-1){0.16}}
\linethickness{0.3mm}
\put(20,30){\line(1,0){30}}
\linethickness{0.3mm}
\put(90,30){\line(0,1){20}}
\linethickness{0.3mm}
\put(90,30){\line(1,0){20}}
\linethickness{0.3mm}
\put(110,30){\line(0,1){20}}
\linethickness{0.3mm}
\put(90,50){\line(1,0){20}}
\put(14,30){\makebox(0,0)[cc]{$s$}}
\put(57,29){\makebox(0,0)[cc]{$y$}}
\put(84,30){\makebox(0,0)[cc]{$z$}}
\put(117,29){\makebox(0,0)[cc]{$y$}}
\put(35,55){\makebox(0,0)[cc]{$z$}}
\put(85,51){\makebox(0,0)[cc]{$s$}}
\put(115,51){\makebox(0,0)[cc]{$r$}}
\end{picture}
\caption{configurations in the undirected graph defining $N$ and $L$ in equation (\ref{bvh77})}\label{yu988}
\end{figure}
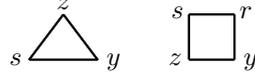

We now follow through the formalism from Section~\ref{secpre} for quantum geodesics, in the graph case. First, from Theorem~\ref{nablachi} and substitution of the forms of $\sigma$ and $\alpha$ in the graph case, we have $\nabla_\cX$ in terms of the $L,N$ as
\begin{align} \label{bvh77}
\nabla_\chi( \chi_{y\leftarrow z} ) &= -\sum_{s:z\to s}   \chi_{y\leftarrow z} \tens\omega_{z\to s} 
+ \sum_{s\to r} \lambda_{s\to r}\, L_{s\to r,z\to y}\,  \chi_{r\leftarrow s} \tens \omega_{s\to z} \cr
& \quad - \sum_s \lambda_{s\to y}\, N_{s\to y,z\to y}\, \chi_{y\leftarrow s} \tens \omega_{s\to z}
\end{align}
where the sums for $L,N$ are over the patterns (possibly degenerate) given by Fig.~\ref{yu988}. It then follows that 
\begin{align} \label{bvh5}
(\id\tens X)\nabla_\chi{}( X) &= -\sum_{s:z\to s} X^{y\leftarrow z}\,   \chi_{y\leftarrow z} \, X^{s\leftarrow z}  \, 
+ \sum_{s\to r}X^{y\leftarrow z}\,  \lambda_{s\to r}\, L_{s\to r,z\to y}\,  \chi_{r\leftarrow s} \,   X^{z\leftarrow s}   \cr
& \quad - \sum_s X^{y\leftarrow z}\,  \lambda_{s\to y}\, N_{s\to y,z\to y}\, \chi_{y\leftarrow s} \,   X^{z\leftarrow s} 
\end{align}
as needed for the geodesic velocity equations.

Next, we consider $\int$ of the form 
\[ \int f=\sum_X \mu(x) f(x)\]
 for all functions $f$ on the vertex set, with `measure' $\mu:X\to \R\setminus\{0\}$ so that $\int$ is hermitian and non-degenerate. We preferably also want $\mu$ to be positive for the usual interpretation of the amplitude flow. 

From the definitions, it is is immediate that a vector field  $X=\sum X^{r\leftarrow s}\,  \chi_{r\leftarrow s}$ is real with respect to $\int$ (i.e. the second half of (\ref{unitarity}) holds) if and only if
\begin{equation}\label{realXgra}
 (X^{y\leftarrow x})^* =-{\mu_y\over \mu_x}\, X^{x\leftarrow y}
\end{equation}
Likewise, it is immediate that
\begin{equation}\label{divintgra}
\mathrm{div}_{\rm \int}(X)(x) = \sum_{y:x\to y} X^{y\leftarrow x} -  \sum_{y:y\to x}\frac{\mu_y}{\mu_x}\  X^{x\leftarrow y} 
\end{equation}
It then follows that 
\begin{align} \label{bvh}
[X,\mathrm{div}_{\int}(X)] &= \sum X^{x\leftarrow y} \, \chi_{x\leftarrow y} \, (\mathrm{div}_{\int}(X)(y)-\mathrm{div}_{\int}(X)(x)) \cr
&=  \sum X^{x\leftarrow y} \, \chi_{x\leftarrow y} \, \Big(
\sum_{z:y\to z} X^{z\leftarrow y} -  \sum_{z:z\to y}\frac{\mu_z}{\mu_y}\  X^{y\leftarrow z}  
- \sum_{z:x\to z} X^{z\leftarrow x} +  \sum_{z:z\to x}\frac{\mu_z}{\mu_x}\  X^{x\leftarrow z} \Big)
\end{align}
as also needed for the geodesic velocity equation for the choice $\kappa={1\over 2}{\rm div}_{\int}(X)$. We put these results together and add the possibility of a driving force $F\in \cX$. 

\begin{proposition}\label{geovelpair} On a graph, the geodesic velocity equation with driving term $F\in \cX$ is
\begin{align*} 
-\dot X^{x\leftarrow y} &= X^{x\leftarrow y} \, \frac12 \, \Big(
-\sum_{z:y\to z} X^{z\leftarrow y} -  \sum_{z:z\to y}\frac{\mu_z}{\mu_y}\  X^{y\leftarrow z}  
- \sum_{z:x\to z} X^{z\leftarrow x} +  \sum_{z:z\to x}\frac{\mu_z}{\mu_x}\  X^{x\leftarrow z} \Big) \cr
 & \quad 
+ \sum_{r,z}X^{r\leftarrow z}\,  \lambda_{y\to x}\, L_{y\to x,z\to r}\,  X^{z\leftarrow y}   
 - \sum_z X^{x\leftarrow z}\,  \lambda_{y\to x}\, N_{y\to x,z\to x}\,  X^{z\leftarrow y} +F^{x\leftarrow y}
\end{align*}
If we suppose that $X$ is real with respect to $\int$ and stays so then then we also simultaneously impose  
\begin{align*} 
-\dot X^{x\leftarrow y} &= X^{x\leftarrow y} \, \frac12 \, \Big(
\sum_{z:x\to z} \frac{\mu_z}{\mu_x}\ X^{x\leftarrow z} +  \sum_{z:z\to x} X^{z\leftarrow x}  
+ \sum_{z:y\to z} \frac{\mu_z}{\mu_y}\  X^{y\leftarrow z} -  \sum_{z:z\to y} X^{z\leftarrow y} \Big) \cr
 & \quad 
- \sum_{r,z} \frac{\mu_r}{\mu_y}\ X^{z\leftarrow r}\,  \lambda_{x\to y}\, L_{z\to r,x\to y}\,  X^{x\leftarrow z}   
 + \sum_z X^{z\leftarrow y}\,  \lambda_{x\to y}\, N_{z\to y,x\to y}\,  X^{x\leftarrow z} -\frac{\mu_x}{\mu_y}\, (F^{y\leftarrow x})^*
\end{align*}
\end{proposition}
\proof We use (\ref{bvh}) and (\ref{bvh5}) to write the 
 geodesic velocity equations in (\ref{veleqX}) as stated. Next we conjugate this equation and suppose that $X$ is real with respect to $\mu$ as in (\ref{realXgra}) to obtain the second equation.  \endproof

The difference between the two can be viewed as an `improved auxiliary equation' that ensures that evolution stays real with respect to $\int$. Setting $F=0$ is the standard approach in \cite{BegMa:cur} but here we make the minimal assumption that $F$ is imaginary but otherwise whatever it needs to be to maintain real evolution for $X$. Finally, given a geodesic velocity field $X$, we have the amplitude flow 
\begin{equation}\label{ampgra}\dot\psi_x=-{1\over 2}\psi_x {\rm div}_{\int}(X)_x-\sum_{p\leftarrow x} (\psi_p-\psi_x) X^{p\leftarrow x}\end{equation}
where ${\rm div}_{\int}$ is the function in (\ref{divintgra}).

The above is the basic framework where $\mu$ is arbitrary. Usually we would take the measure adapted to the geometry and the natural way to do this is to ask further that  $\int$ is divergence compatible. Then we have a further structure of a $*$ operation on vector fields. 

\begin{proposition}\label{propHvv} 
\[ {\rm div}(\chi_{p\leftarrow q})=\delta_q- \delta_p\,  \frac{ \lambda_{p\to q }  }  {  \lambda_{q\to p}  }  \sum_{ z \, : \,  p\to z  } L_{   p\to q  , z\to p  }
\, \lambda_{p\to z} \]
Moreover, $\int$ is divergence compatible if
\begin{align*} 
 \sum_{ z \, : \,  p\to z  } L_{   p\to q  , z\to p  }
\, \lambda_{p\to z} =\frac{\mu_q}{\mu_p}\  \frac  {  \lambda_{q\to p}  }  { \lambda_{p\to q }  } 
\end{align*}
In this case the canonical $*$ on $\cX$ is
\[
(\chi_{p\leftarrow q})^* = -  \frac{ \lambda_{p\to q }  }  {  \lambda_{q\to p}  }  \sum_{ z \, : \,  p\to z  } 
L_{  z\to p, p\to q   }
\, \lambda_{p\to z} \, \chi_{q\leftarrow p}.
\]
\end{proposition}
\proof Starting with
\[
\mathfrak{g}_2(\omega_{w\to z})=-\lambda_{z\to w}\,\chi_{w\leftarrow z} 
\]
we use (\ref{sigmachiinv}) to find
\begin{align*}
\sigma_\cX{}^{-1}(\chi_{p\leftarrow q}\tens\omega_{q\to r}) &= - (\id\tens \mathfrak{g}_2) \,\sigma_{\Omega^1}(\omega_{p\to q}\tens\omega_{q\to r})\, \lambda_{q\to p}{}^{-1} \\
&=   - (\id\tens \mathfrak{g}_2) \, \sum_{ z \, : \,  p\to z\to r  } \lambda_{p\to q }\,L_{   p\to q  , z\to r  }\,  \omega_{p\to z} \tens \omega_{z\to r}
\, \lambda_{q\to p}{}^{-1} \\
&=    \sum_{ z \, : \,  p\to z\to r  } \lambda_{p\to q }\,L_{   p\to q  , z\to r  }\,  \omega_{p\to z} \tens \chi_{z\leftarrow r}
\, \lambda_{q\to p}{}^{-1}\,\lambda_{r\to z} 
\end{align*}
and so $\tilde\ev:=\ev\circ\sigma_\cX^{-1}$ is 
\begin{align*}
\tilde\ev(\chi_{p\leftarrow q}\tens\omega_{q\to p}) 
&= \delta_p\,  \frac{ \lambda_{p\to q }  }  {  \lambda_{q\to p}  }  \sum_{ z \, : \,  p\to z  } L_{   p\to q  , z\to p  }
\, \lambda_{p\to z} 
\end{align*}
where $\chi_{y\leftarrow x}$ is the dual basis to $\omega_{x\to y}$. Next, we observe that any connection on an inner calculus, we have 
\[
\mathrm{div}(X)=\ev(\theta\tens X)-\tilde{\ev}(X\tens\theta)-(\tilde{\ev}\tens\ev)(e_i\tens\alpha(e^i)\tens X)
\]
for  dual bases, which in our case means  $\coev(1)=\sum_{r,s}\chi_{r\leftarrow s}\tens \omega_{s\to r}$. This gives the formula for ${\rm div}(\chi_{p\leftarrow q})$ stated on noting that the $\alpha$ term does not contribute. From this, the divergence-compatibility immediately follows. 

Finally, for  $*$ in a basis vector, we compute
\begin{align*}
\ev(\omega_{r\to s}\tens (\chi_{p\leftarrow q})^*) &= \tilde\ev(\chi_{p\leftarrow q}\tens (\omega_{r\to s})^*)^*
= -\tilde\ev(\chi_{p\leftarrow q}\tens \omega_{s\to r})^*\cr
&= -\delta_{s,q}\,\delta_{r,p}\, \Big(\delta_p\,  \frac{ \lambda_{p\to q }  }  {  \lambda_{q\to p}  }  \sum_{ z \, : \,  p\to z  } L_{   p\to q  , z\to p  }
\, \lambda_{p\to z} \Big)^*
= -\delta_{s,q}\,\delta_{r,p}\, \delta_p\,  \frac{ \lambda_{p\to q }  }  {  \lambda_{q\to p}  }  \sum_{ z \, : \,  p\to z  } 
L_{  z\to p, p\to q   }
\, \lambda_{p\to z} 
\end{align*}
which implies the result stated. \endproof

In our case of $\varsigma=\id$ it is shown in \cite{BegMa:cur} that $X=X^*$ is the same as $X$ being real with respect to $\int$, so we are solving the same pair of equations as in Proposition~\ref{geovelpair} but have more structure underlying the reality condition. In the other hand, this case can be quite restrictive.

\section{Quantum geodesics on the 4-star graph}\label{secstar}

We now apply the theory above to the $n$-star, proceeding from the pair of equations in Proposition~\ref{geovelpair}. In fact it is sufficient to set $x=0$ (setting $y=0$ gives the conjugate equation which we impose anyway), so we have to solve the velocity equation
 \begin{align*} 
-\dot X^{0\leftarrow y} &= X^{0\leftarrow y} \, \frac12 \, \Big(
- X^{0\leftarrow y} - \frac{\mu_0}{\mu_y}\  X^{y\leftarrow 0}  
- \sum_{i} X^{i\leftarrow 0} +  \sum_{i}\frac{\mu_i}{\mu_0}\  X^{0\leftarrow i} \Big) \cr
 & \quad 
+ \sum_{r}X^{r\leftarrow 0}\,  \lambda_{y\to 0}\, L_{y\to 0,0\to r}\,  X^{0\leftarrow y}   +F^{0\leftarrow y}\ ,\cr
\end{align*} and, for $X$ to remain real with respect to $\int$, at the same time
\begin{align*}
-\dot X^{0\leftarrow y}  &= X^{0\leftarrow y} \, \frac12 \, \Big(
\sum_{i} \frac{\mu_i}{\mu_0}\ X^{0\leftarrow i} +  \sum_{i} X^{i\leftarrow 0}  
+ \frac{\mu_0}{\mu_y}\  X^{y\leftarrow 0} -  X^{0\leftarrow y} \Big) \cr
 & \quad 
- \sum_{z} \frac{\mu_0}{\mu_y}\ X^{z\leftarrow 0}\,  \lambda_{0\to y}\, L_{z\to 0,0\to y}\,  X^{0\leftarrow z}    -\frac{\mu_0}{\mu_y}\, (F^{y\leftarrow 0})^*.
\end{align*}
The difference of these is the `improved auxiliary condition' 
 \begin{align*} 
 & X^{0\leftarrow y} \, \frac12 \, \Big(
- X^{0\leftarrow y} - \frac{\mu_0}{\mu_y}\  X^{y\leftarrow 0}  
- \sum_{i} X^{i\leftarrow 0} +  \sum_{i}\frac{\mu_i}{\mu_0}\  X^{0\leftarrow i} \Big) \cr
 & \quad 
+ \sum_{r}X^{r\leftarrow 0}\,  \lambda_{y\to 0}\, L_{y\to 0,0\to r}\,  X^{0\leftarrow y}   +F^{0\leftarrow y}  \cr
&= X^{0\leftarrow y} \, \frac12 \, \Big(
\sum_{i} \frac{\mu_i}{\mu_0}\ X^{0\leftarrow i} +  \sum_{i} X^{i\leftarrow 0}  
+ \frac{\mu_0}{\mu_y}\  X^{y\leftarrow 0} -  X^{0\leftarrow y} \Big) \cr
 & \quad 
- \sum_{z} \frac{\mu_0}{\mu_y}\ X^{z\leftarrow 0}\,  \lambda_{0\to y}\, L_{z\to 0,0\to y}\,  X^{0\leftarrow z}    -\frac{\mu_0}{\mu_y}\, (F^{y\leftarrow 0})^*
\end{align*}
which simplifies to
 \begin{align*} 
 F^{0\leftarrow y}+ \frac{\mu_0}{\mu_y}\, (F^{y\leftarrow 0})^* 
&= X^{0\leftarrow y} \,  \Big(
  \sum_{i} X^{i\leftarrow 0}  
+ \frac{\mu_0}{\mu_y}\  X^{y\leftarrow 0}  \Big) \cr
 & \quad 
- \sum_{i} \frac{\mu_0}{\mu_y}\ X^{i\leftarrow 0}\,  \lambda_{0\to y}\, L_{i\to 0,0\to y}\,  X^{0\leftarrow i}    
- \sum_{i}X^{i\leftarrow 0}\,  \lambda_{y\to 0}\, L_{y\to 0,0\to i}\,  X^{0\leftarrow y}   
\end{align*}

We now suppose that $n=4$ and that $\nabla$ is the QLC as found in Theorem~\ref{thmsolstar}. Here $s=-1$ and the auxiliary equation becomes
 \begin{align*} 
 F^{0\leftarrow y}+ \frac{\mu_0}{\mu_y}\, (F^{y\leftarrow 0})^* &= X^{0\leftarrow y} \,  \Big(
  \sum_{i} X^{i\leftarrow 0}  
+ \frac{\mu_0}{\mu_y}\  X^{y\leftarrow 0}  \Big)  \quad 
- {1\over 2} \sum_{i} \frac{\mu_0}{\mu_y}\ X^{i\leftarrow 0}\,   X^{0\leftarrow i}    
- {1\over 2} \sum_{i}X^{i\leftarrow 0}\, \frac{ \lambda_{y\to 0}} { \lambda_{0\to i}  }   \,   X^{0\leftarrow y}   \cr
 & \quad 
+ \frac{\mu_0}{\mu_y}\ X^{y\leftarrow 0}\,   X^{0\leftarrow y}    
+  X^{y\leftarrow 0}\,  \frac{ \lambda_{y\to 0}} { \lambda_{0\to y}  }    X^{0\leftarrow y}   \\
 &= X^{0\leftarrow y} \,  
  \sum_{i} X^{i\leftarrow 0}  \Big( 1  -  \frac{ \lambda_{0\to y}} { \lambda_{0\to i}  }    \Big)
-  \frac12\,  \frac{\mu_0}{\mu_y}  \sum_{i}X^{i\leftarrow 0} \,  X^{0\leftarrow i}    
+ 2\, \Big( 1    + \frac{\mu_0}{\mu_y}\Big)
\ X^{0\leftarrow y} \,   X^{y\leftarrow 0} 
\end{align*}
If we assume $F$ is imaginary with respect to $\int$ then this amounts to requiring 
 \begin{align} \label{Fstar}
 F^{0\leftarrow y}=  \frac{\mu_0}{\mu_y}\, (F^{y\leftarrow 0})^*
&=  \frac12\,  X^{0\leftarrow y} \,  
  \sum_{i} X^{i\leftarrow 0}  \Big( 1  -  \frac{ \lambda_{0\to y}} { \lambda_{0\to i}  }    \Big)
-  \frac14\,  \frac{\mu_0}{\mu_y}  \sum_{i}X^{i\leftarrow 0} \,  X^{0\leftarrow i}    
+ \Big( 1    + \frac{\mu_0}{\mu_y}\Big)
\ X^{0\leftarrow y} \,   X^{y\leftarrow 0}.
\end{align}
We interpret $F$ as an external force defined by this and needed so as to keep $X$ real with respect to $\int$ during the evolution. 

\begin{proposition} With the driving term (\ref{Fstar}), the geodesic velocity equation on the 4-star becomes 
\begin{align*}-\dot X^{0\leftarrow y}&= \frac12 \, X^{0\leftarrow y} \, \bigg( 
- X^{0\leftarrow y}  
+  \sum_{i}\frac{\mu_i}{\mu_0}\  X^{0\leftarrow i}  
-\sum_{i} (X^{0\leftarrow i})^*\,\frac{\mu_i\, \lambda_{0\to y}} {\mu_0 \lambda_{0\to i}  } 
+ \Big( 2 \frac{\mu_y}{\mu_0} - 1\Big) \,   (X^{0\leftarrow y})^* \bigg)
+  \frac14\,  \sum_{i}  \frac{\mu_i}{\mu_y} \, | X^{0\leftarrow i} |^2  
\end{align*}
for four complex fields $X^{0\leftarrow i}$.
\end{proposition}
\proof We now put the found $F$ into our original geodesic velocity equation to give
\begin{align*} 
-\dot X^{0\leftarrow y} &= X^{0\leftarrow y} \, \frac12 \, \Big(
- X^{0\leftarrow y} - \frac{\mu_0}{\mu_y}\  X^{y\leftarrow 0}  
- \sum_{i} X^{i\leftarrow 0} +  \sum_{i}\frac{\mu_i}{\mu_0}\  X^{0\leftarrow i} \Big) \cr
 & \quad 
+ \sum_{r}X^{r\leftarrow 0}\,  \lambda_{y\to 0}\, L_{y\to 0,0\to r}\,  X^{0\leftarrow y}   +F^{0\leftarrow y}   \cr
&= X^{0\leftarrow y} \, \bigg( \frac12 \, \Big(
- X^{0\leftarrow y} - \frac{\mu_0}{\mu_y}\  X^{y\leftarrow 0}  
- \sum_{i} X^{i\leftarrow 0} +  \sum_{i}\frac{\mu_i}{\mu_0}\  X^{0\leftarrow i} \Big) \cr
 & \quad 
+ \sum_{i}X^{i\leftarrow 0}\,  \lambda_{y\to 0}\, L_{y\to 0,0\to i}  +    \frac12 \,  
  \sum_{i} X^{i\leftarrow 0}  \Big( 1  -  \frac{ \lambda_{0\to y}} { \lambda_{0\to i}  }    \Big)  
+ \Big( 1    + \frac{\mu_0}{\mu_y}\Big) \,   X^{y\leftarrow 0} \bigg)
-  \frac14\,  \frac{\mu_0}{\mu_y}  \sum_{i}X^{i\leftarrow 0} \,  X^{0\leftarrow i}   \cr
&= \frac12 \, X^{0\leftarrow y} \, \bigg( 
- X^{0\leftarrow y}  
+  \sum_{i}\frac{\mu_i}{\mu_0}\  X^{0\leftarrow i}  
+\sum_{i}X^{i\leftarrow 0}\, \Big( 2\,  \lambda_{y\to 0}\, L_{y\to 0,0\to i}  
-    \frac{ \lambda_{0\to y}} { \lambda_{0\to i}  }   \Big)
+ \Big( 2    + \frac{\mu_0}{\mu_y}\Big) \,   X^{y\leftarrow 0} \bigg)
\\ &\quad  -  \frac14\,  \frac{\mu_0}{\mu_y}  \sum_{i}X^{i\leftarrow 0} \,  X^{0\leftarrow i}   \cr
&= \frac12 \, X^{0\leftarrow y} \, \bigg( 
- X^{0\leftarrow y}  
+  \sum_{i}\frac{\mu_i}{\mu_0}\  X^{0\leftarrow i}  
+\sum_{i}X^{i\leftarrow 0}\,\frac{ \lambda_{0\to y}} { \lambda_{0\to i}  } \, \Big( 2\,  \frac{\lambda_{y\to 0}}{\lambda_{0\to y}}
\, \lambda_{0\to i} \,L_{y\to 0,0\to i}  
-    1  \Big)
+ \Big( 2    + \frac{\mu_0}{\mu_y}\Big) \,   X^{y\leftarrow 0} \bigg)
\\ &\quad  -  \frac14\,  \frac{\mu_0}{\mu_y}  \sum_{i}X^{i\leftarrow 0} \,  X^{0\leftarrow i}   \cr
&= \frac12 \, X^{0\leftarrow y} \, \bigg( 
- X^{0\leftarrow y}  
+  \sum_{i}\frac{\mu_i}{\mu_0}\  X^{0\leftarrow i}  
+\sum_{i}X^{i\leftarrow 0}\,\frac{ \lambda_{0\to y}} { \lambda_{0\to i}  } \, \Big( 4
\, \lambda_{0\to i} \,L_{y\to 0,0\to i}  
-    1  \Big)
+ \Big( 2    + \frac{\mu_0}{\mu_y}\Big) \,   X^{y\leftarrow 0} \bigg)
 \\ & \quad -  \frac14\,  \frac{\mu_0}{\mu_y}  \sum_{i}X^{i\leftarrow 0} \,  X^{0\leftarrow i}   \cr
&= \frac12 \, X^{0\leftarrow y} \, \bigg( 
- X^{0\leftarrow y}  
+  \sum_{i}\frac{\mu_i}{\mu_0}\  X^{0\leftarrow i}  
+\sum_{i}X^{i\leftarrow 0}\,\frac{ \lambda_{0\to y}} { \lambda_{0\to i}  } 
+ \Big(  \frac{\mu_0}{\mu_y} - 2\Big) \,   X^{y\leftarrow 0} \bigg)
 -  \frac14\,  \sum_{i}  \frac{\mu_0}{\mu_y} \, X^{i\leftarrow 0} \,  X^{0\leftarrow i} 
\end{align*}
for the value of $L$, which we write as shown. \endproof

\begin{example}\label{ex4star}\rm Figure~\ref{fig4star} shows solutions for the geodesic velocity equation on the 4-star graph for constant $\lambda_{0\to i}$ (independent of  $i$), constant $\mu$  and initial $X^{0\leftarrow i}(0)=\delta_{i,1}$. A special feature of $n=4$ (not true for $n=2,3$) is that both the velocity and amplitude coefficients remain real numbers if they start real, so we stick to this for simplicity. For the constant measure case we then have $X^{i\leftarrow 0}=-X^{0\leftarrow i}$ at all $s$ for reality with respect to $\int$. The velocity equation for the stated initial conditions then reduces to 
\[ X^{0\leftarrow 1}=\xi,\quad X^{0\leftarrow 2}=X^{0\leftarrow 3}=X^{0\leftarrow 3}=\xi-1;\quad \dot\xi=- \xi^2 + {3\over 2} \xi -{3\over 4},\quad \xi(0)=1,\]  
which is solved by
\[\xi(s):=\frac{1}{4} \left(\sqrt{3} \tan \left(\frac{\pi }{6}-\frac{\sqrt{3} }{4}s\right)+3\right),\quad s\in [0,{8\pi\over 3\sqrt{3}}),\]
where the $\tan$ first blows up at $s={8\pi\over 3\sqrt{3}}$. We also compute
\[ {\rm div}_{\int}(X)_0=6-8 \xi ,\quad {\rm div}_{\int}(X)_1=2\xi,\quad {\rm div}_{\int}(X)_2={\rm div}_{\int}(X)_3={\rm div}_{\int}(X)_4=2(\xi-1),\]
which gives the amplitude flow equations (\ref{ampgra}) as
\[ \dot \psi_0=\xi \psi_1+ 3 (\xi-1)\psi_i,\quad \dot\psi_1=-\xi \psi_0,\quad \dot\psi_i=-(\xi-1)\psi_0,\quad i=2,3,4.\]
The numerical solution for this is plotted for initial $\psi_1(0)=1$ and zero elsewhere. Here $|\psi|^2_x$ is the probability to find a particle at vertex $x$ and we see that this shifts from outer vertex $x=1$ to the central vertex $x=0$ and then equally to the other outer vertices $x=2,3,4$. It then continues to shift back and forth in an oscillatory fashion between 0 and the outer vertices in unison and with divergent frequency as we approach the singularity at $s={8\pi\over 3\sqrt{3}}$. One can check that $\sum_x |\psi_x|^2=1$ for all $s$ to within machine precision. 

Moreover, since this is a unitary evolution, there is necessarily an effective time-dependent Hamiltonian, namely we can write
\[ \dot \psi=\imath H_s\psi,\quad H_s=-\imath\begin{pmatrix} 0 & \xi & 3(\xi-1)\\ -\xi & 0 & 0\\ 1-\xi & 0 &0\end{pmatrix}\]
(taking $\psi$ with components $\psi_0,\psi_1,\psi_2$). Here at each $s$, $H_s$ has one real zero mode and two complex modes, 
\[ \begin{pmatrix}0 \\ 3(1-\xi)\\ \xi\end{pmatrix},\quad \begin{pmatrix}\pm \imath\lambda \\  \xi \\ \xi-1\end{pmatrix};\quad \lambda={\sqrt{3}\over \sqrt{2}\sqrt{1+\sin \left(\frac{\pi}{6}+ {\sqrt{3}s
\over 2}\right)}}\ge {\sqrt{3}\over 2}\]
with eigenvalues $0,\mp \lambda$ respectively, with $\lambda$ blowing up as $s$ approaches the end of its range. Our plotted solution starts off real at  $s=0$ as the sum of the two complex modes and then remains real,  but a general solution would be complex. \begin{figure}
\[ \includegraphics[scale=0.85]{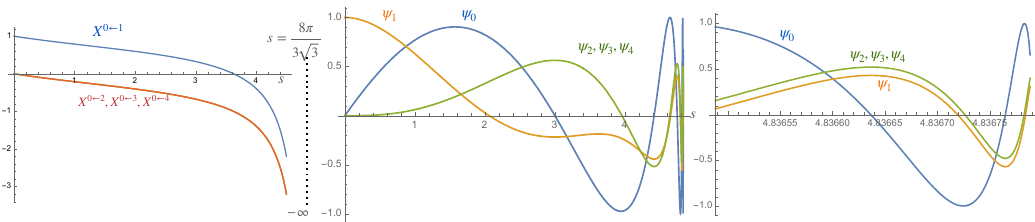}\]
\[\includegraphics[scale=1.1]{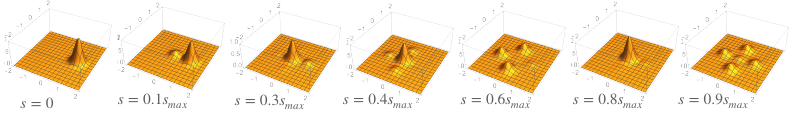}\]
\caption{\label{fig4star} Example~\ref{ex4star} of quantum geodesic flow on the 4-star graph.  The geodesic velocity field has initial value $X^{0\leftarrow 1}=1$ only,  and  its components blow up in a finite time $s_{max}={8\pi\over 3\sqrt{3}}$. We plot the amplitude flow  $\psi_x$ where $|\psi_x|^2$ is the probability to find the particle at vertex $x$. The 3D visualisation puts the central vertex $0$ at $(0,0)$, the outer vertices at $(\pm 1,0),(0,\pm 1)$ and Gaussian interpolates the values of $|\psi|^2$. We see the particle starting at outer vertex $(1,0)$, moving to $(0,0)$ and then moving to the three other outer vertices. It then increasingly rapidly oscillates between the external nodes in unison and the centre. }
\end{figure}
\end{example}

\begin{remark}\rm In the above, we have not tried to choose $\int$ to be divergence compatible and instead let $\mu$ (and $\lambda_{x\to y}$) be arbitrary. If we require divergence compatibility then for the 4-star this amounts to 
\[ \sum_j \lambda_{0\to j}=0,\quad \mu_i=-{1\over 2} \mu_0\]
for the measure to exist and its resulting form. This means that the $\lambda$ cannot all be positive, i.e. the metric cannot have `Euclidean signature' everywhere. And we see that the measure is not positive either, which is a problem for the probabilistic interpretation. Aside from these issues, the system itself can be similarly solved to the example. \end{remark}

\section{Quantum geodesics on Cayley graphs}\label{seccayley}

In this section,  we look at special graphs where the vertices are elements $x\in G$ of a group  (so $A=\C(G)$) and arrows are of the form $x\to ya$ where $a\in \CC$ a finite set of generators. In this case we define
\[ e^a=\sum_{x\to xa} \omega_{x\to xa},\quad  e^a f=R_a(f)e^a,\quad \extd f=\sum_a(\del_a f)e^a,\quad \del_a=R_a-\id,\quad R_a(f)(x)=f(xa).\]
The vector space spanned by the $\{e^a\}$ is $\Lambda^1$ and is a left basis for $\Omega^1=A.\Lambda^1$ as a free module. We assume throughout that $\CC$ is closed under inversion. Then we  have an inner $*$-calculus with
\[\theta=\sum_a e^a,\quad e^a{}^*=-e^{a^{-1}}.\]

If $\CC$ is ad-stable then the calculus is bicovariant, and in this case there is a canonical $\Omega_{wor}(G)$ defined by a braiding $\Psi(e^a\tens e^b)=e^{aba^{-1}}\tens e^a$. This tends to be a natural quotient of $\Omega_{min}$ defined by the graph structure alone. The choice of $\Omega^2$ is relevant to conditions for torsion of QLC or to define a WQLC.

\subsection{QRG on Cayley graphs revisited} \label{seccay}

A Hermitian metric on $\Omega^1$ has the form 
\[ \<e^a,\overline{e^b}\>=h_{ab}=\delta_{a,b}h_a,\]
where the natural condition is for $h^{ab}$ to be a positive Hermitian matrix, so here $h_a>0$, as functions on the group. The diagonal form comes from the observation that corresponding corresponding regular metric $g$ has to be central to have a bimodule inverse. Because of the minus sign in $e^a{}^*$, this is  
\[ g=-g_a e^a\tens e^{a^{-1}},\quad  (e^a,e^b)=-\delta_{a,b^{-1}}h_a,\quad   h_a={1\over R_a(g_{a^{-1}})}. \]
The edge symmetric case has $R_a(g_{a^{-1}})=g_a$ but we do not necessarily assume this. If we are only interested in the round brackets version, we usually omit the minus signs here as in \cite[Chap.~1]{BegMa}, since the QLC does not depend on the overall normalisation of the metric.

Now consider a left connection on $\Omega^1$ of the form
\begin{align} \label{nablacay}
 \nabla e^a=-\sum_{b,c}\Gamma^a{}_{bc}e^b\tens e^c.
\end{align}
We set matrices $\Gamma$ and $\Xi$ with entries in $\Omega^1$ to be
\begin{equation}\label{Xi}
\Gamma_{ab} =  \sum_c \Gamma^a{}_{cb}e^c\ ,\quad \Xi_{ab} = \sum_c ( \Gamma^a{}_{cb}+\delta_{a,b})e^c\ ,\quad
(\Xi_{ab})_c=\Xi^a{}_{cb}:= \Gamma^a{}_{cb}+\delta_{a,b}
\end{equation}
so we have $\Xi=\Gamma+\theta I$ regarded as a 1-form valued matrix, where $I$ is the identity matrix. 

\begin{lemma} \label{hermgroup}  
A connection on a Cayley graph as in (\ref{nablacay}) is hermitian-metric compatible if and only if  matrix $\Xi h$ is antihermitian where we matrix multiply $h_{ab}$. Given its diagonal form, the condition is
\[
 h_b R_{c} (\Xi^a{}_{c^{-1}b}{}^*)
=    \Xi^b{}_{ca} R_c( h_a).
\]
\end{lemma}
\proof We use the formula from \cite[Def.~8.33]{BegMa}
\[
\extd h = -\Gamma h - h \Gamma^*
\]
and the fact that the inner element $\theta$ is antihermitian. In the case where $h$ is diagonal, this implies
\[
( \Xi_{ab}h_b)^*  = \sum_c ( \Xi^a{}_{cb}e^c h_b)^*= - \sum_c h_b e^{c^{-1}} \Xi^a{}_{cb}{}^*
= - \sum_c h_b e^{c} \Xi^a{}_{c^{-1}b}{}^*
=-\Xi_{ba}h_a=  - \sum_c  \Xi^b{}_{ca}e^c h_a
\]
which is the explicit equation stated on moving $e^c$ to the right in the last expressions. Similarly in the converse direction.
 \endproof
 
 Next, we ask when $\nabla$ is a bimodule connection. This was answered under the assumption that $\nabla$ as left invariant so that $\Gamma^a{}_{bc}$ are $\C$-valued in \cite[Prop.~3.75]{BegMa} but we show now that this assumption is not necessary for the same conclusion.
 
 \begin{proposition}\label{propbimod} On a Cayley graph, $\nabla$ given by (\ref{Xi}) is a bimodule connection if and only if 
\[
a^{-1}bc \notin \mathcal{C}\cup\{e\} \implies  \Xi^a{}_{bc}=0.
\]
In this case,
\[
\sigma(e^a\tens e^b) =  \sum_{c,d\in \mathcal{C}:cd=ab} \Xi^a{}_{cd}\,  e^c\tens e^d,\quad \alpha(e^a)=0\]
 in the decomposition (\ref{nablainner}), i.e. a bimodule connection on a Cayley graph is necessarily inner.
\end{proposition}
\proof (1) We calculate  
\begin{align*}
\nabla(e^a\,\delta_s)-(\nabla e^a)\delta_s &= \nabla(\delta_{sa^{-1}} e^a) + \sum_{b,c}\Gamma^a{}_{bc}\delta_{sc^{-1}b^{-1}}e^b\tens e^c \\
&= \sum_{c\in\mathcal{C}} ( \delta_{sa^{-1} c^{-1}}  - \delta_{sa^{-1}}   )   e^c\tens e^a +\sum_{b,c\in\mathcal{C}}  (\delta_{sc^{-1}b^{-1}} - \delta_{sa^{-1}} ) \Gamma^a{}_{bc}e^b\tens e^c
\end{align*}
and for this to equal $\sigma(e^a\tens\extd \delta_s)$,  we require
\begin{align*}
\sum_{b\in\mathcal{C}} (\delta_{sb^{-1}a^{-1}}-\delta_{sa^{-1}}) \sigma(e^a\tens e^b)
&= \sum_{c\in\mathcal{C}} ( \delta_{sa^{-1} c^{-1}}  - \delta_{sa^{-1}}   )   e^c\tens e^a +\sum_{b,c\in\mathcal{C}}  (\delta_{sc^{-1}b^{-1}} - \delta_{sa^{-1}} ) \Gamma^a{}_{bc}e^b\tens e^c\ .
\end{align*}
Now taking $\delta_{sa^{-1}}$ times this and summing over $s\in G$ gives
\begin{align} \label{sumone}
\sum_{b\in\mathcal{C}} \sigma(e^a\tens e^b)
&= \sum_{c\in\mathcal{C}}  e^c\tens e^a +\sum_{b,c\in\mathcal{C}: bc\neq a}   \Gamma^a{}_{bc}e^b\tens e^c\ .
\end{align}
from which we derive the requirement
\begin{align} \label{sumtwo}
\sum_{b\in\mathcal{C}} \delta_{sb^{-1}a^{-1}}\sigma(e^a\tens e^b)
&= \sum_{c\in\mathcal{C}}  \delta_{sa^{-1} c^{-1}}  e^c\tens e^a +\sum_{d,c\in\mathcal{C} : dc\neq a} \delta_{sc^{-1}d^{-1}} \Gamma^a{}_{dc}e^d\tens e^c   \ .
\end{align}
Now taking $\delta_{sb^{-1}a^{-1}}$ times this and summing over $s\in G$ gives 
\begin{equation}\label{sigmacay} \sigma(e^a\tens e^b) = \sum_{c,d\in \mathcal{C}:cd=ab} \Gamma^a{}_{cd} e^c\tens e^d +
\left\{\begin{array}{cc}e^{aba^{-1}}\tens e^a & aba^{-1}\in\mathcal{C} \\0 & aba^{-1}\notin\mathcal{C}\end{array}\right.
=  \sum_{c,d\in \mathcal{C}:cd=ab} (\Gamma^a{}_{cd}+\delta_{a,d})  e^c\tens e^d\end{equation}
as stated, and this formula is supposed to give a bimodule map. Taking $\delta_{sx^{-1}}$ times (\ref{sumtwo}) for $a^{-1}x\notin\mathcal{C}$ and summing over $s\in G$ gives
\begin{align} \label{sumthree}
0
&= \sum_{d\in\mathcal{C}}  \delta_{da,x}  e^d\tens e^a +\sum_{d,c\in\mathcal{C} : dc\neq a} \delta_{dc,x} \Gamma^a{}_{dc}e^d\tens e^c   \ .
\end{align}
For $x=a$ this just gives $0=0$, and for the remaining case $a^{-1}x\notin\mathcal{C}\cup\{e\}$ we get
\begin{align} \label{sumfour}
0
&= \sum_{d\in\mathcal{C}}  \delta_{da,x}  e^d\tens e^a +\sum_{d,c\in\mathcal{C} } \delta_{dc,x} \Gamma^a{}_{dc}e^d\tens e^c   \ .
\end{align}
so for a given $d\in\mathcal{C}$ we get
\begin{align} \label{sumfive}
0
&=  \delta_{da,x}  e^a +\sum_{c\in\mathcal{C} } \delta_{dc,x} \Gamma^a{}_{dc}e^c   \ .
\end{align}
which gives $\delta_{dc,x}( \Gamma^a{}_{dc}+\delta_{c,a})=0$, which is the condition in the statement. 
A brief check shows that (\ref{sumone}) is then satisfied and that we obtain a bimodule connection by the formula stated.   

(2) Next, we compute
\begin{align*}
\nabla e^a- \theta\tens e^a+\sigma(e^a\tens\theta)&=-\sum_{c,d}\Gamma^a{}_{cd}e^c\tens e^d-\sum_b e^b\tens e^a+\sum_{b=a^{-1}cd}\Big(\Gamma^a{}_{cd}+\delta_{a,d}\Big)e^c\tens e^d\\
&=\sum_{a^{-1}cd\notin\CC}\delta_{a,d}e^c\tens e^d-\sum_b e^b\tens e^a+\sum_{a^{-1}cd\in \CC}\delta_{ad}e^c\tens e^d=0\end{align*}
where we used (\ref{sigmacay}) for the first expression. The first sum with $\Gamma$ is unconstrained but the other sum is constrained, so the difference is a sum over $a^{-1}cd\notin\CC$. But for such values we use (\ref{nablacay}). 
  \endproof

That $\alpha=0$ is implicit (but explicitly stated) in the treatment of left-invariant bimodule connections on Cayley graphs in \cite{BegMa} but we see that this too holds in general. 
 
 \begin{lemma} \label{urptt} A bimodule  connection on a Cayley graph is  $*$-compatible if and only if
 \[
\sum_{c,d\in \mathcal{C}:cd=ab} R_{ (ab)^{-1}} ( \Xi^a{}_{cd} {}^*) 
\Xi^{d^{-1}}{}_{rs} = \delta_{b^{-1},r}\, \delta_{a^{-1},s}.
\]
and is automatically $*$-preserving. \end{lemma}
\proof Using our formula from Proposition~\ref{propbimod}, 
\[
\dag \sigma(e^a\tens e^b) = \dag \sum_{c,d\in \mathcal{C}:cd=ab}  \Xi^a{}_{cd}   e^c\tens e^d = 
\sum_{c,d\in \mathcal{C}:cd=ab}  e^{d^{-1}} \tens  e^{c^{-1}}  \Xi^a{}_{cd} {}^* = 
\sum_{c,d\in \mathcal{C}:cd=ab} R_{ (ab)^{-1}} ( \Xi^a{}_{cd} {}^*) e^{d^{-1}} \tens  e^{c^{-1}}
\]
and then the  $*$-compatibility condition is 
\[
\sigma\dag \sigma(e^a\tens e^b)= \sum_{r,s\in \mathcal{C}:rs=(ab)^{-1}} \sum_{c,d\in \mathcal{C}:cd=ab} R_{ (ab)^{-1}} ( \Xi^a{}_{cd} {}^*) 
\Xi^{d^{-1}}{}_{rs} e^r\tens e^s = e^{b^{-1}} \tens e^{a^{-1}}
\]
so for all $a,b\in\mathcal{C}$ we have the condition stated. Also, because $\alpha=0$, being $*$-compatible implies $*$-preserving by \cite[Thm.~8.11]{BegMa}. 
\endproof

Similarly,  since $\alpha=0$, for $\nabla$ to be torsion free is equivalent to being torsion compatible\cite[Thm.~8.11]{BegMa}, i.e. the condition $\wedge(\id+\sigma)=0$, or explicitly
\[ \sum_{cd=ab}\Xi^a{}_{cd}\, e^c\wedge e^d+ e^a\wedge e^b=0\]
for all $a,b\in \CC$. The wedge product here depends on $\Omega^2$, e.g for the Woronowicz exterior algebra when the calculus is bicovariant, it means the same expression with $\tens$ in place of $\wedge$  needs to be invariant under $\Psi$. For $G=S_3$, there is a 1-parameter family of left-invariant torsion free $*$-preserving WQLCs\cite[Prop.~6.12]{BegMa}.

Finally, if $\nabla$ is given by (\ref{nablacay}), we will need the associated right connection $\nabla_\cX$ on vector fields. Here we let $\{f_a\}$ be a dual (right) basis to the (left) basis $\{e^a\}$ in the sense $f_b(e^a)=\ev(e^a\tens f_b)=\delta^a{}_b$. A general vector field here is $X=  f_a X^a \in \cX$. Then 
\[ \nabla_\cX f_a=\sum_{b,c} f_b\tens\Gamma^b{}_{ca}e^c\] 
with braiding  $\sigma_\cX:\Omega^1\tens_A \cX\to  \cX \tens_A  \Omega^1$   given by
\[
\sigma_\cX=(\id\tens\id\tens\ev)(\id\tens\sigma\tens\id)(\sum_c f_c\tens e^c\tens\id\tens\id)
\]
resulting in
\[
\sigma_\cX(e^a\tens f_b) = \sum_{r,c \, : \, rb=ca}  R_{c^{-1}} (\Xi^c{}_{rb} )f_c\tens e^r\ .
\]

\subsection{The integral state and divergence}

As for any graph we define a state $\int f=\sum_{x} \mu_x f(x)$ for some non-vanishing and preferably positive weight function $\mu:G\to \R\setminus\{0\}$. In this case, the associated divergence is computed from
\begin{align*}  \int(X(\extd f)) &= \sum_{a\in\mathcal{C}}  \sum_{x\in G} X((R_{a}(f) - f )  e^a)(x)\mu(x) \cr
&=  \sum_{a\in\mathcal{C}}  \sum_{x\in G} \big(  X^a(x)\mu(x) R_{a}(f)(x)  -  X^a(x)\mu(x) f(x)\big) \cr
&=  \sum_{a\in\mathcal{C}}  \sum_{x\in G} \big(R_{a^{-1}}(X^a\mu)(x)  - X^a(x)\mu(x) \big)f(x) \cr
&= \int\Big( \sum_{a\in\mathcal{C}}  \big(  R_{a^{-1}}(X^a\mu) \mu^{-1}  - X^a \big) f \Big)
\end{align*}
so that 
\begin{equation}\label{divint}
\mathrm{div}_{\int} X =  \sum_{a\in\mathcal{C}}  \big(  X^a- R_{a^{-1}}(X^a\mu) \mu^{-1}   \big),
\end{equation}
which on a basis vector is just $\mathrm{div}_{\int} (f_d )=  1- R_{d^{-1}}(\mu) \mu^{-1}$. Moreover, a vector field $X$ is 
defined as real with respect to $\int$ if the second half of the `unitarity' condition (\ref{unitarity}), which 
  amounts to 
\[
\int(X(\delta_x e^a)^*)=\int(X(e^a{}^*\delta_x))=-\int(X(e^{a^{-1}}\delta_x))=-\int(X(\delta_{x a}e^{a^{-1}}))
\]
and gives  $\int(\delta_x X^a{}^*)=-\int(\delta_{x a}X^{a^{-1}})$ or
\begin{equation}\label{XrealG} X^a(x)^*=-X^{a^{-1}}(xa){\mu(xa)\over \mu(x)}\end{equation}
for all $a\in \CC, x\in G$ as the condition to be real with respect to $\int$. In this case, $\kappa={1\over 2}{\rm div}_{\int}(X)$ is real-valued and ensures the other half of (\ref{unitarity}). 

We can then proceed to quantum geodesics on groups for  any measure, but as in \cite{BegMa:cur} a natural choice is to ask for $\int$ to be divergence compatible in the sense $\int {\rm div}(X)=0$ for all $X\in \cX$, or equivalently that 
${\rm div}_{\int}={\rm div}$.

\begin{proposition}  \label{urp} If $\nabla$ is metric-compatible then the geometric divergence of a left vector field $X$ is 
\[ {\rm div}(X)=\sum_{a} X^{a^{-1}}-\sum_{a,b}g_a R_a(X^{a^{-1}}) {\Xi^a{}_{b,b^{-1}}\over R_b(g_{b^{-1}}  )  }\]
and $\int$ is divergence compatible  if and only if 
\[ {R_a(\mu)\over\mu}=g_a \sum_b {\Xi^a{}_{b,b^{-1}}\over R_b(g_{b^{-1}})}.\]
In this case, the induced $*$-operation on $\cX$ is
\[ (X^*)^a=- R_a(X^{a^{-1}})^* g_a \sum_b {\Xi^a{}_{b,b^{-1}}{}^*\over R_b(g_{a^{-1}})}= -R_a(X^{a^{-1}})^* {R_a(\mu)\over\mu}. \]
\end{proposition}
\proof The corresponding 1-form to the left vector field $X$ via metric is 
\[ \mathfrak{g}_2^{-1}(X)=g^1X(g^2)=-\sum_{a,b}g_a e^a\ (f_b X^b)(e^{a^{-1}})=-\sum_a g_a e^a X^{a^{-1}}=- \sum_ a g_a R_a(X^{a^{-1}}) e^{a}.\]
Then by (\ref{div2}), 
\begin{align*}{\rm div}(X)=(\ ,\ )\nabla \mathfrak{g}_2^{-1}(X)&=\sum_a (\ ,\ )\left(-\sum_b\del_b (g_a R_{a}(X^{a^{-1}}))e^b\tens e^a+ g_a R_a(X^{a^{-1}})\sum_{b,c}\Gamma^a{}_{bc}e^b\tens e^c\right)\\
&=\sum_a {\del_{a^{-1}} (g_a R_a(X^{a^{-1}}))\over R_{a^{-1}}(g_a)}- \sum_{a,b} g_a R_a(X^{a^{-1}}){\Gamma^a{}_{bb^{-1}}\over R_b(g_{b^{-1}})}\\
&=\sum_a X^{a^{-1}}- \sum_a {g_a R_a(X^{a^{-1}})\over R_{a^{-1}}(g_a)}-\sum_{a,b}  g_a R_a(X^{a^{-1}}) {\Gamma^a{}_{bb^{-1}}\over R_b(g_{b^{-1}})}
\end{align*} 
which we write as stated. Then
\[ \int {\rm div}(X)=\sum_G \sum_a R^a(X^{a-1})  \left( R_a(\mu)-\sum_{b}\mu g_a  {\Xi^a{}_{b,b^{-1}}\over R_b(g_{b^{-1}}}\right)\]
since we could apply $R_a$ to the first term given the $\sum_G$. For this to vanish for all $X$, we need the condition stated. 

Finally, when $\mu$ obeys this condition then by Theorem~\ref{thmXstar}, we have a $*$ operation on $\cX$. This is characterised by
\begin{align*} (X^*)^a&=\ev(e^a\tens X^*)=[\ev\circ\sigma_\cX^{-1}(X\tens (e^a)^*)]^*=-[(\ ,\ )\sigma(\mathfrak{g}_2(X)\tens e^{a^{-1}})]^*\\
&=[\sum_{cd=ba} g_b R_b(X^{b^{-1}})\Xi^b{}_{cd} (e^c,e^d)]^*\end{align*}
on using the second half of (\ref{sigmachiinv}) in the theorem, followed by $\mathfrak{g}_2(X)$ as above and $\sigma$ from Proposition~\ref{propbimod}. We then use the value of $(e^c,e^d)$ to obtain the first stated expression, which we then recognise as the ratio $R_a(\mu)/\mu$. 
\endproof

We see that the condition $X^*=X$ is the same as (\ref{XrealG}) for reality with respect to $\int$, but is now arising as the self-adjoint elements with respect to an involution on $\cX$. 

\subsection{Geodesic velocity equation on a discrete group}

We assume that we are in the setting where we have fixed $\int$ compatible with the geometric divergence and hence have  defined $*$ for a real vector field.  We now consider the equations (\ref{veleqX}) and their $*$ in the discrete group case. 

\begin{proposition}\label{veleqCay} On a Cayley graph, the geodesic velocity equation for a time-dependent vector field $X$ and with driving term $F\in\cX$ is
\[  -\dot X^a = -X^a\del_a(\kappa) +  \sum_{b,d}  \Xi^a{}_{bd} R_b(X^d)X^b-  \sum_{b} X^a X^b+F^a\]
Moreover, for $X$ to remain real with respect to $\int$ and $\kappa={1\over 2}{\rm div}_{\int}(X)$ (assuming that the connection preserves the hermitian metric) under the assumption that 
$F$ is imaginary with respect to $\int$ then we require
\begin{align*}
-2\,F^a&=  -X^a\, \sum_{b} \Big( R_{b^{-1}}(X^b\mu)/\mu + R_a(    X^b)  \Big) +  \sum_{b,d}    R_a(h_{a^{-1}} ) \, R_{ab^{-1}}(    \Xi^{d^{-1}}{}_{ba^{-1}} \, R_{d^{-1}}(X^{d}\mu)
\, X^{b}    /h_{d^{-1}})   / \mu
 \cr
&\quad +  \sum_{b,d}  \Xi^{a}{}_{bd} R_b(X^d)X^b  
\end{align*}
The velocity equation is
\begin{align*}
-2\dot X^a &=X^a \sum_b \Big(-X^b+R_{ab^{-1}}(X^b\mu)/R_a(\mu)\Big) +  \sum_{b,d}  \Xi^{a}{}_{bd} R_b(X^d)X^b \cr
&   -  \sum_{b,d}    R_a(h_{a^{-1}} ) \, R_{ab^{-1}}(    \Xi^{d^{-1}}{}_{ba^{-1}} \, R_{d^{-1}}(X^{d}\mu)
\, X^{b}    /h_{d^{-1}})   / \mu
\end{align*}
\end{proposition}
\proof For the first part,
\[
[X,\kappa]=f_a X^a\kappa-\kappa\, f_a X^a =f_a X^a(\kappa-R_a(\kappa))
\]
\begin{align*}
\nabla_\cX(X) &= \nabla_\cX(f_d X^d)=f_d\tens\extd X^d +  \sum_{b,a} f_a\tens \Gamma^a{}_{bd} e^bX^d
= \sum_{d,b} f_d\tens (R_b(X^d)-X^d)e^b+  \sum_{b,a,d} f_a\tens \Gamma^a{}_{bd} R_b(X^d)e^b \cr
&= \sum_{b,a,d} f_a\tens \Xi^a{}_{bd} R_b(X^d)e^b-  \sum_{a,b} f_a\tens X^a e^b
\end{align*}
so
\[
(\id\tens X)\nabla(X)=  \sum_{b,a,d} f_a \Xi^a{}_{bd} R_b(X^d)X^b-  \sum_{a,b} f_aX^a X^b
\]
and (\ref{veleqX}) becomes, on addition of a forcing term
\[
\dot X^a + X^a(\kappa-R_a(\kappa)) +  \sum_{b,d}  \Xi^a{}_{bd} R_b(X^d)X^b-  \sum_{b} X^a X^b+F^a=0
\]
as stated. For the second part, we define 
\[
B^a= X^a\, (\kappa-R_a(\kappa))+  \sum_{b,d}  \Xi^{a}{}_{bd} R_b(X^d)X^b-  \sum_{b} X^a X^b
\]
so that our geodesic velocity equation is $-\dot X^a =B^a+F^a$. We let $X^a$ be is real with respect to $\int$ and $F^a$ is imaginary (i.e.\ $\mathrm{i}\, F^a$ is real with respect to $\int$). Conjugating the original equation and suitably subtracting gives us
\[ F^a=-{1\over 2}(B^a+{R_a((B^{a^{-1}})^*\,\mu)\over \mu}).\]
Next, for our particular $B^a$, we have 
\[
R_a(B^{a^{-1}}\,\mu)=R_a(X^{a^{-1}}\mu)\, R_a(\kappa-R_{a^{-1}}(\kappa)) +  \sum_{b,d}  R_a(\mu\,\Xi^{a^{-1}}{}_{bd} R_b(X^d)X^b) 
-  \sum_{b} R_a( \mu\,X^{a^{-1}} X^b)
\]
and conjugating and using the reality of $X$ (and therefore $\kappa$) gives
\[
R_a(B^{a^{-1}}\,\mu)^*/\mu= X^a\, (\kappa-R_a(\kappa))+  \sum_{b,d}  R_a(\mu\,\Xi^{a^{-1}}{}_{bd} R_b(X^d)X^b)^*/\mu
+  \sum_{b}  X^a\, R_a(  X^b)^*
\]
Now
\begin{align} \label{ghu}
B^a+R_a(B^{a^{-1}}\,\mu)^*/\mu &= 2X^a\, (\kappa-R_a(\kappa)) +  \sum_{b,d}  R_a(\mu\,\Xi^{a^{-1}}{}_{bd} R_b(X^d)X^b)^*/\mu
+  \sum_{b}  X^a\, R_a(  X^b)^*   \cr
&\quad 
+  \sum_{b,d}  \Xi^{a}{}_{bd} R_b(X^d)X^b-  \sum_{b} X^a X^b
\end{align}
The terms not containing $\Xi$ in (\ref{ghu}) are
\begin{align} \label{ghu2}
X^a\, (2\kappa  -  \sum_{b}  X^b -R_a(2\kappa
-  \sum_{b}     X^b{}^*)  )
\end{align}
and using the divergence (\ref{divint}) and reality condition (\ref{XrealG}) this becomes
\begin{align} \label{ghu3}
X^a\, (\sum_{b}  X^b{}^* -R_a(\sum_{b}     X^b)  ) = -X^a\, \sum_{b} \Big( R_{b^{-1}}(X^b\mu)/\mu + R_a(    X^b) 
\Big)
\end{align}
Using the condition for hermitian metric preservation, the two terms containing $\Xi$ in (\ref{ghu}) are
\begin{align*}
& \sum_{b,d}    R_a(h_{a^{-1}} ) \, R_{ab}(\Xi^d{}_{b^{-1}a^{-1}})\,    R_a(\mu\, R_b(X^d)X^b)^*/( \mu\, R_{ab}(h_d))
+  \sum_{b,d}  \Xi^{a}{}_{bd} R_b(X^d)X^b \cr
&= - \sum_{b,d}    R_a(h_{a^{-1}} ) \, R_{ab}(\Xi^d{}_{b^{-1}a^{-1}})\,    R_a(R_b(X^d
{}^*\, X^{b^{-1}}\mu)   )/( \mu\, R_{ab}(h_d))
+  \sum_{b,d}  \Xi^{a}{}_{bd} R_b(X^d)X^b  \cr
&=  \sum_{b,d}    R_a(h_{a^{-1}} ) \, R_{ab}(\Xi^d{}_{b^{-1}a^{-1}})\,    R_{ab}(R_d(X^{d^{-1}}\mu)
\, X^{b^{-1}})   /( \mu\, R_{ab}(h_d))
+  \sum_{b,d}  \Xi^{a}{}_{bd} R_b(X^d)X^b  \cr
&=  \sum_{b,d}    R_a(h_{a^{-1}} ) \, R_{ab}(    \Xi^d{}_{b^{-1}a^{-1}} \, R_d(X^{d^{-1}}\mu)
\, X^{b^{-1}}    /h_d)   / \mu
+  \sum_{b,d}  \Xi^{a}{}_{bd} R_b(X^d)X^b  \cr
&=  \sum_{b,d}    R_a(h_{a^{-1}} ) \, R_{ab^{-1}}(    \Xi^{d^{-1}}{}_{ba^{-1}} \, R_{d^{-1}}(X^{d}\mu)
\, X^{b}    /h_{d^{-1}})   / \mu
+  \sum_{b,d}  \Xi^{a}{}_{bd} R_b(X^d)X^b  
\end{align*}
and combining these gives the stated answer.
\endproof 

We can also write the final result in Proposition~\ref{veleqCay} as 
\begin{align}
-2\dot X^a &=X^a \sum_b \Big(-X^b+R_{ab^{-1}}(X^b\mu)/R_a(\mu)\Big) +  \sum_{b,d}  \Xi^{a}{}_{bd} R_b(X^d)X^b  \cr
&   -  \sum_{b,d}     R_{a}(    \Xi^{a^{-1}}{}_{b^{-1} d^{-1}} {}^*) \,  R_{ab^{-1}}(   R_{d^{-1}}(X^{d}\mu)
\, X^{b}  )   / \mu \label{veleqCaystar}
\end{align}
by rewinding the last step of the proof where we used hermitian metric compatibility. In the proposition, we added a natural driving force $F$ needed to simultaneously enforce both the velocity equation and its conjugate for $X$ real with respect to $\int$.  The measure $\mu$ is also arbitrary but, as before, a natural choice is for the integral to be divergence compatible so that ${\rm div}_{\int}={\rm div}$, the geometric divergence. 
 
 After solving for the geodesic velocity flow $X^a$, we then solve for the `amplitude flow' on time-dependent $\psi\in \C(G)$, which is therefore
\[ \dot \psi= -\sum_a(\del_a \psi)X^a- \psi\kappa\]
with probability density $\rho=\bar \psi\psi$ when suitably normalised. It follows from the theory that
\[ \int \dot \rho=0\]
so that $\rho$ can be normalised to a probability measure with respect to $\int$. 

\section{Quantum geodesic scattering on the integer lattice line}\label{secZ}

Here, we see how the Cayley graph theory applies to $G=\Z$, the integer line graph with $\CC=\{\pm 1\}$. The simplest case of the theory was recently studied in \cite{BegMa:Z} in the edge-symmetric case and it was found that divergence-compatibility forces the the curvature to be zero, i.e. the metric weights must form a geometric sequence. Our theory above is more general and allows us to cover generic metrics with functions $g_\pm(i)=g_{i\to i\pm 1}$, but we stick to the 
edge symmetric case where $g_-=R_-(g_+)$ so that there is one functions $g_i:=g_+(i)=g_{i\to i+1}=g_{i+1\to i}$ say for the metric weight attached to each edge. We let
\[ \rho_\pm=R_\pm({g_\pm\over g_\mp}),\quad \rho_+={R_+(g)\over g}=:\rho,\quad \rho_-=R_-^2(\rho^{-1})\]
and use the unique $*$-preserving QLC\cite{ArgMa1}, 
\[ \nabla e^\pm=(1-\rho_\pm) e^\pm\tens e^\pm,\quad \sigma(e^\pm\tens e^\pm)=\rho_\pm e^\pm\tens e^\pm ,\quad \sigma(e^\pm\tens e^\mp)=e^\mp\tens e^\pm.  \]
This also descends to $\Z_n$, which has a unique $*$-preserving QLC when $n\ne 4$. When $n=4$ with the smaller $\Omega_{wor}$, there are some further such $*$-preserving QLCs (it is not known if there is a unique QLC for $n=4$ if we use the bigger calculus $\Omega_{min}$). 

Here the Christoffel symbols in the basis $e^\pm$ have only basis have only two non-zero entries $\Gamma^\pm{}_{\pm\pm}=\rho_\pm-1$, hence $\nabla_\cX$ in the dual basis $f_\pm$ is 
\[ \nabla_\cX f_\pm = f_\pm\tens (\rho_\pm-1) e^\pm,\quad \sigma_\cX(e^\pm\tens f_\pm)= f_\pm\tens \rho_\pm e^\pm,\quad \sigma_\cX(e^\pm\tens f_\mp)= f_\mp\tens e^\pm.\]
One can check that these are indeed left and right connections using the commutation rules
\[ e_\pm f= R_\pm(f) e^\pm,\quad  f f_\pm= f_\pm R_\pm(f)\]
and metric compatible as claimed. The structure constants from the point of view of the theory in Section~\ref{seccayley} are
\begin{equation}\label{XiZ} \Xi^{\pm}{}_{\pm\pm}=\rho_\pm, \quad \Xi^{\pm}{}_{\mp\pm}=1\end{equation}
 with all others zero. 

We now proceed in the general setting without assuming divergence compatibility, so both the edge-symmetric metric and the measure are arbitrary. From (\ref{divint})-(\ref{XrealG}), we have
\[ {\rm div}_{\int}(X)= -{1\over\mu}(\del_-(\mu X^+)+\del_+(\mu X^-)), \quad X^\pm{}^*=-{R_\pm(\mu X^\mp)\over\mu}=-\mu_\pm R_\pm X^\mp,\quad \mu_\pm={R_\pm(\mu)\over\mu}\]
for the divergence and reality with respect to $\int$.  Setting $\kappa={1\over 2}{\rm div}_{\int}(X)$ we have
\[2 \del_-(\kappa)=-\del_-({1\over\mu}\del_-(\mu X^+))+{1\over R_-(\mu)}\del_-(\mu X^-)+{1\over\mu}\del_+(\mu X^-).\]
The geodesic velocity equation as in Proposition~\ref{veleqCay} can be computed as 
\[ \dot X^\pm = \del_\pm(\kappa)X^\pm + (1-\rho_\pm)X^\pm X^\pm -\rho_\pm\del_\pm(X^\pm)X^\pm- \del_\mp(X^\pm)X^\mp-F^\pm\]
where we add a driving term. We assume this is imaginary with respect to $\int$ and apply * to the $\dot X^+$ equation, divide by $-\mu_+$ and then apply  $R_-$, to obtain 
\[ \dot X^-=-(\del_-\kappa)X^--(1-R_-(\rho))R_-(\mu_+)X^-{}^2+R_-(\rho)\del_+(R_-(\mu_+)X^-)X^-+{1\over R_-(\mu_+ R_-(\mu_+))}\del_-(R_-(\mu_+)X^-)R_-^2(X^+)+F^-\]
which we subtract from the $\dot X^-$ equation to obtain
\begin{align*} 2F^-&=2\del_-(\kappa)X^-+(1-\rho_-)X^-{}^2+{1-R_-(\rho)\over\mu_-}X^-{}^2 - \rho_-(\del_- X^-)X^--R_-(\rho)\del_+({X^-\over\mu_-})X^-\\
&\quad -\del_+(X^-)X^+-\mu_-R_-(\mu_-)\del_-({X^-\over\mu_-})R_-^2(X^+)\\
&={1\over R_-(\mu)}\del_-(\mu X^-)+{1\over\mu}\del_+(\mu X^-)
+(1-\rho_-)X^-{}^2+{1-R_-(\rho)\over\mu_-}X^-{}^2 - \rho_-(\del_- X^-)X^--R_-(\rho)\del_+({X^-\over\mu_-})X^-\\
&\quad -\del_-({1\over\mu}\del_-(\mu X^+)) -\del_+(X^-)X^+-\mu_-R_-(\mu_-)\del_-({X^-\over\mu_-})R_-^2(X^+)\\
&=(1-\rho_-)R_-(X^-)X^-+ \mu_+(1-R_-(\rho))(R_+ X^-)X^-- (\mu_-R_--\id)\del_-(X^+R_+(X^-))
\end{align*}
after a lot of cancellations. We used $R_-(\mu_+)=1/\mu_-$. We then apply * to obtain the other half of
\[ 2 F^\pm=(1-\rho_\pm)R_\pm(X^\pm)X^\pm+ \mu_\mp(1-R_\pm(\rho_\mp))(R_\mp X^\pm)X^\pm- (\mu_\pm R_\pm-\id)\del_\pm (X^\mp R_\mp (X^\pm))\]
We now put this driving force term into the velocity equation to obtain finally\begin{align*}\dot X^\pm &= \del_\pm(\kappa)X^\pm + (1-\rho_\pm)X^\pm X^\pm -\rho_\pm\del_\pm(X^\pm)X^\pm- \del_\mp(X^\pm)X^\mp\\
&\quad -{1\over 2}((1-\rho_\pm)R_\pm(X^\pm)X^\pm+ \mu_\mp(1-R_\pm(\rho_\mp))(R_\mp X^\pm)X^\pm- (\mu_\pm R_\pm-\id)\del_\pm (X^\mp R_\mp (X^\pm)))\end{align*}
which simplifies to
\begin{equation}\label{veleqZ} 2\dot X^\pm=\left(\mu_\mp R_\pm(\rho_\mp)R_\mp- \rho_\pm R_\pm+ (1-{1\over\mu_\pm})\id\right)(X^\pm)X^\pm + \left(\mu_\pm R_\pm- R_\pm(\mu_\pm)\id\right)(X^\pm)R_\pm^2(X^\mp)-(\del_\mp X^\pm)X^\mp.\end{equation}
We only have to solve one of these as, by construction, the other is the conjugate. We have given a direct derivation but one can verify that the same results are obtained from Proposition~\ref{veleqCay} or more easily from (\ref{veleqCaystar}). Since $\Xi$ in our case is real, this comes down to
\begin{align*}-2\dot X^+ &=X^+ \sum_b \Big(-X^b+R_{1-b}(X^b\mu)/R_+(\mu)\Big) +  \sum_{b}  \Xi^{+}{}_{b+} R_b(X^+)X^b \cr
&   -  \sum_{b}     R_{+}(    \Xi^{-}{}_{b^{-1} -} ) \,  R_{1-b}(   R_{-}(X^{+}\mu)
\, X^{b}  )   / \mu \end{align*}
which combined with (\ref{XiZ})  give the same as (\ref{veleqZ}). After solving for $X^+$, say, with any initial distribution on $\Z$,  we are then free to choose  any initial wave function $\psi$ and solve the amplitude flow equations
\[  \dot \psi = -X_s(\extd \psi )-\psi \kappa_s= -(\del_+\psi)X^+- (\del_-\psi)X^--\psi\kappa. \]

\begin{figure}
\[\includegraphics[scale=.75]{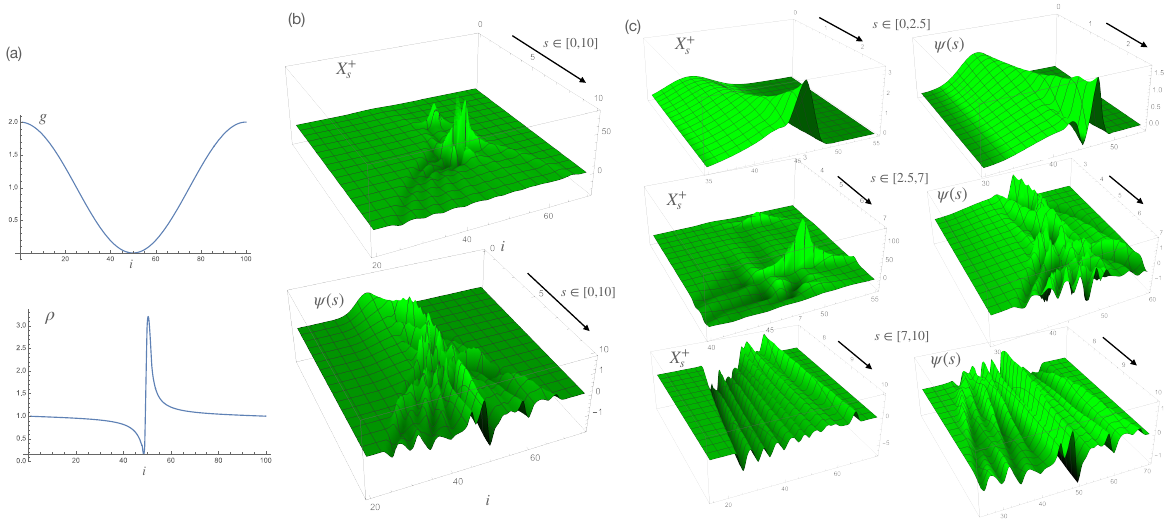}\]
\caption{(a) Metric with a minimum at $i=50$. (b) Resulting velocity field $X^+$ and wave function $\psi$ (shown interpolated) as functions of geodesic time $s$, for initial Gaussians centred  $i=40$. Part (c) shows the same in close-up as three segments of geodesic time. \label{figZcos}}
\end{figure}

\begin{example}\rm  So far, the measure and metric are arbitrary. We now exhibit an explicit solutions and for simplicity we take  $\mu$ constant so that $\mu_\pm=1$. Then the geodesic velocity equation for general metric becomes
\[ 2\dot X^+=({1\over R_-(\rho)}R_- - \rho R_+)(X^+)X^++ \del_+(X^+) R_+^2(X^-)- \del_-(X^+)X^-,\quad X^-=-R_-(X^+{}^*).\]
For our example, we will stick with $X^+$ real-valued (then it stays real during the evolution) in which case
\[ X^-=-R_-(X^+), \quad \kappa=-{1\over 2}(\del_- X^++\del_+ X^-)=-\del_- X^+.\]
Our equation in terms of $X^+$ alone becomes
\begin{equation}\label{veleqZm1} \dot X^+={1\over 2}\left(({1\over R_-(\rho)}-1)R_- +(1- \rho) R_+\right)(X^+)X^++ {1\over 2}(R_-(X^+)^2-  R_+(X^+)^2).\end{equation}

 For the quantum metric we choose one that dips as a cosine (sampled on the integer lattice) as shown in Fig.~\ref{figZcos}(a) along with the corresponding function $\rho$. We then solve the velocity equations (\ref{veleqZm1}) with initial value of $X^+$ at time $s=0$ a Gaussian located at $i=40$ at the start of dip. This moves towards the centre then becomes large and turns into a wave packet as shown in  Fig.~\ref{figZcos}(b)-(c).  The values are interpolated for purposes of visualisation. 
 
 Next, the amplitude equation for real-valued $X^+$ becomes
\[ \dot\psi= -(\del_+ \psi)X^+ + (\del_-\psi)R_-(X^+)+ \psi \del_-(X^+)=-R_+(\psi)X^++ R_-(\psi)R_-(X^+).\]
which we solve for the previously found $X^+$, with interpolated results  shown in Fig.~\ref{figZcos}(b)-(c) for an initial Gaussian centred at $i=40$ for maximum effect. (If the initial $\psi$ Gaussian were to be centred to the left or the right of this then would keep its shape for the main part but with similar ripples appearing where part of it overlaps the evolving $X^+$ distribution. \end{example}

\begin{example}\rm  An initial analysis of quantum-geodesics on $\Z$ in \cite{BegMa:Z} shows that divergence compatibility of $\int f=\sum_i f\mu$ is  possible if and only if $\rho$  is constant so the metric is exponential in $i$ (the sequence of metric values is a geometric sequence) and in which case the natural measure is $\mu=g$. From the QRG point of view, this is exactly the case where the Riemann curvature vanishes. This case leads by Theorem~\ref{thmXstar}  to
\[ X^+{}^*=-\rho R_+(X^-),\quad X^-{}^*=-{1\over\rho}R_-(X^+),\quad {\rm div}(X)=(\id- {R_-\over \rho})X^++ (\id-\rho R_+)X^-\]
and the geodesic velocity equation (\ref{veleqZ}) and $F$ become
\begin{equation}\label{veleqZrho}  2\dot X^+=\big( {R_-\over \rho^2}-\rho R_++ (1-{1\over \rho})\big)(X^+)X^+ + \rho\del_+(X^+)R^2_+(X^-)-\del_-(X^+)X^-\end{equation}
\[ 2F^+=(1-\rho) \big(R_+- {R_-\over\rho^2}\big)(X^+)X^+- (R_+-{\id\over\rho})\del_-(|X^+|^2)\]
where we used that  $\mu_\pm=\rho^{\pm 1}$ and
\[ X^- R_- (X^+)= R_-(X^+ R_+ (X^-))=-{R_-(|X^+|^2)\over\rho}.\]
to interpret the second term of $F$. 

\begin{figure}
\[\includegraphics[scale=.85]{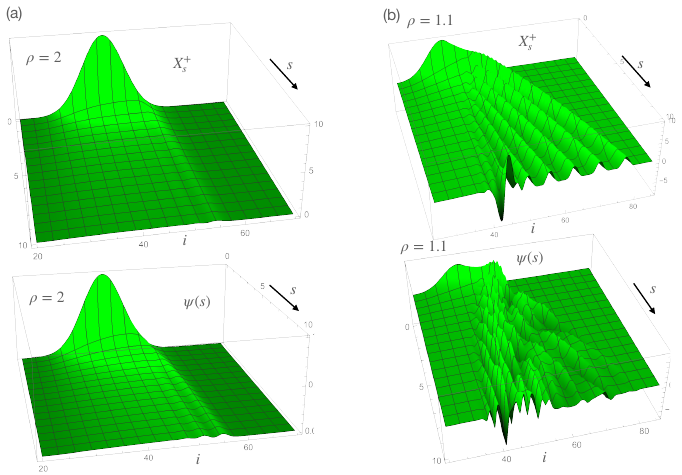}\]
\caption{Quantum geodesics for an exponentially increasing metric $g_i=g_0\rho^i$ with (a) $\rho=2$ and (b) $\rho=1.1$. In both cases, the geodesic velocity $X^+$ and the wave function $\psi$ (shown interpolated) are  initially Gaussians centred at $i=40$.  \label{figZflat}}
\end{figure}

If we assume for simplicity that $X^+$ is real-valued then we can use $X^-=-{R_-(X^+)\over\rho}$ to write (\ref{veleqZrho}) entirely in terms of $X^+$, to obtain
\[ \dot X^+={1-\rho\over 2} \left( {R_-\over\rho^2}+ R_+ - {\id\over\rho}\right)(X^+)X^++{1\over 2}\left( {R_- (X^+)^2\over\rho}-R_+ (X^+)^2\right)\]
and $\kappa=(\id- {R_-\over\rho})X^+$. The amplitude flow similarly in the real-valued case  then becomes  
\[  \dot \psi= -R_+(\psi)X^++ R_-(\psi){R_-(X^+)\over\rho}.\]
Some numerical solutions where $X^+$ starts at $s=0$ as a Gaussian at $i=40$ and $\psi$ starts (for maximum effect) also starts as a Gaussian at $i=40$ are shown in Fig.~\ref{figZflat}. Interestingly, for a metric that grow less rapidly so that $\rho$ is close to 1, we see large oscillations appearing in the evolution. The unitarity is with respect to the measure $\mu_i=\rho^i$ so it is $\sum_{i\in \Z} \rho^i |\psi_i|^2$ which is preserved under the evolution (as one can check to within numerical accuracy).

The results in this example are complementary to \cite{BegMa:Z}, where we imposed $F=0$ (as an auxiliarly condition needed for geodesic flow without imaginary driving forces). In that case it is shown, at least for $\rho=1$, that real-valued $X^+$ are not possible aside from some trivial cases (but rather it should have  constant absolute value). This then leads to the phenomenon that initially real $\psi$ become complex wave functions during the geodesic evolution.

\end{example}

 \end{document}